\documentclass[12pt,leqno]{amsart}
\usepackage{amssymb,amsthm,verbatim,enumerate,hyperref,ifthen}
\usepackage[mathscr]{eucal}
\usepackage[utf8]{inputenc}
\usepackage[T1]{fontenc}
\usepackage{graphics}
\usepackage{graphicx}
\def\N{\mathbb{N}}
\def\R{\mathbb{R}}

\def\Z{\mathbb{Z}}

\def\sup{\mathop{\mbox{\rm sup}}}
\def\dist{\mathop{\mbox{\rm dist}}}

\def\meas{\mathop{\mbox{\rm meas}}}

\def\epsilon{\varepsilon}
\def\phi{\varphi}
\def\eps{\varepsilon}

\def\/{|\!|\!|}
\newtheorem{theorem}{Theorem}

\newtheorem*{theorem*}{Theorem}
\def\Thm#1#2{\ifthenelse{\equal{#1}{*}}{\begin{theorem*}#2\end{theorem*}}
  {\begin{theorem}\label{T#1}#2\end{theorem}}}
\newtheorem{Atheorem}{Theorem}

\def\THM#1#2{\begin{Atheorem}\label{T#1}#2\end{Atheorem}}
\def\thm#1{Theorem~\ref{T#1}}
\newtheorem{proposition}[theorem]{Proposition}
\newtheorem*{proposition*}{Proposition}
\def\Prp#1#2{\ifthenelse{\equal{#1}{*}}{\begin{proposition*}#2\end{proposition*}}
             {\begin{proposition}\label{P#1}#2\end{proposition}}}
\def\prp#1{Proposition~\ref{P#1}}
\newtheorem{corollary}[theorem]{Corollary}
\newtheorem*{corollary*}{Corollary}
\def\Cor#1#2{\ifthenelse{\equal{#1}{*}}{\begin{corollary*}#2\end{corollary*}}
             {\begin{corollary}\label{C#1}#2\end{corollary}}}

\newtheorem{lemma}[theorem]{Lemma}
\newtheorem*{lemma*}{Lemma}
\def\Lem#1#2{\ifthenelse{\equal{#1}{*}}{\begin{lemma*}#2\end{lemma*}}
             {\begin{lemma}\label{L#1}#2\end{lemma}}}
\def\lem#1{Lemma~\ref{L#1}}
\newtheorem{example}[theorem]{Example}
\newtheorem*{example*}{Example}
\long\def\Exa#1#2{\ifthenelse{\equal{#1}{*}}{\begin{example*}\rm #2\end{example*}}
            {\begin{example}\label{Ex#1}\rm #2\end{example}}}

\newtheorem{problem}[subsection]{Problem}

\theoremstyle{definition}
\newtheorem{definition}[theorem]{Definition}
\newtheorem*{definition*}{Definition}
\def\Defi#1#2{\ifthenelse{\equal{#1}{*}}{\begin{definition*}#2\end{definition*}}
      {\begin{definition}\label{D#1}#2\end{definition}}}

\newtheorem{remark}[theorem]{Remark}
\newtheorem*{remark*}{Remark}
\def\Rem#1#2{\ifthenelse{\equal{#1}{*}}{\begin{remark*}#2\end{remark*}}
             {\begin{remark}\label{R#1}#2\end{remark}}}

\def\eq#1{{\rm(\ref{E#1})}}
\def\Eq#1#2{\ifthenelse{\equal{#1}{*}}
  {\begin{equation*}\begin{aligned}#2\end{aligned}\end{equation*}}
  {\begin{equation}\begin{aligned}\label{E#1}#2\end{aligned}\end{equation}}}

\begin{document}
\begin{flushright}
Cent. Eur. J. Math. \textbf{10}(3) (2012), 1017–1041. \\
\href{http://dx.doi.org/10.2478/s11533-012-0027-5}{doi: 10.2478/s11533-012-0027-5} \\[1cm]
\end{flushright}

\title[Approximate Hermite--Hadamard inequality]
{Implications between approximate convexity properties  
and approximate Hermite--Hadamard inequalities}
\author[J.\ Mak\'o]{Judit Mak\'o}
\author[Zs. P\'ales]{Zsolt P\'ales}
\address{Institute of Mathematics, University of Debrecen,
H-4010 Debrecen, Pf.\ 12, Hungary}
\email{\{makoj,pales\}@science.unideb.hu} \subjclass[2000]{Primary 39B22, 39B12} 
\keywords{convexity, approximate convexity, lower and upper Hermite--Hadamard inequalities}
\thanks{This research has been supported by the Hungarian
Scientific Research Fund (OTKA) Grant NK81402 and by the
TÁMOP-4.2.1/B-09/1/KONV-2010-0007, TÁMOP-4.2.2/B-10/1-2010-0024 projects.
These projects are co-financed by the European Union and the European Social Fund.}

\begin{abstract}
In this paper, the connection between the functional inequalities
$$
  f\Big(\frac{x+y}{2}\Big)\leq\frac{f(x)+f(y)}{2}+\alpha_J(x-y) \qquad (x,y\in D)
$$
and
$$
\int_0^1f\big(tx+(1-t)y\big)\rho(t)dt
\leq\lambda f(x)+(1-\lambda)f(y) +\alpha_H(x-y) \qquad (x,y\in D)
$$
is investigated,
where $D$ is a convex subset of a linear space, $f:D\to\R$, 
$\alpha_H,\alpha_J:D-D\to\R$ are even functions, $\lambda\in[0,1]$,
and $\rho:[0,1]\to\R_+$ is an integrable nonnegative function with $\int_0^1\rho(t)dt=1$.
\end{abstract}

\maketitle

\section{Introduction}

Throughout this paper $\R$, $\R_+$, $\N$ and $\Z$ denote the sets of real, nonnegative real, natural
and integer numbers, respectively. Let $X$ be a real linear space and $D\subset X$ be a convex set.
Denote by $D^*$ the difference set of $D$:
\Eq{*}{
  D^*:=D-D:=\{x-y\mid x,y\in D\}.
}
Of course, $D^*$ is convex and $0\in D^*$.
It is well-known (see \cite{Had1893}, \cite{MitLac85}, \cite{Kuc85}, and
\cite{NicPer06}, \cite{DraPea00}) that convex functions $f:D\to\R$
satisfy the so-called lower and upper Hermite--Hadamard inequalities
\Eq{00l}{
  f\Big(\frac{x+y}{2}\Big)\leq\int_0^1 f\big(tx+(1-t)y\big)dt
  \qquad(x,y\in D),
}
and
\Eq{00}{
  \int_0^1 f\big(tx+(1-t)y\big)dt\leq\frac{f(x)+f(y)}{2}
  \qquad(x,y\in D),
}
respectively. The converse is also known to be true (cf.\ \cite{NicPer03},
\cite{NicPer06}), i.e., if a function $f:D\to\R$ which is continuous over the 
segments of $D$ satisfies \eq{00l} or \eq{00}, then it is also convex.

More generally, it is easy to see that the $\varepsilon$-convexity of $f$
(cf.\ \cite{HyeUla52}), i.e., the validity of
\Eq{*}{
  f\big(tx+(1-t)y)\leq tf(x)+(1-t)f(y)+\varepsilon
  \qquad(x,y\in D,\,t\in[0,1]),
}
implies the following $\varepsilon$-Hermite--Hadamard inequalities
\Eq{nrsl}{
     f\Big(\frac{x+y}{2}\Big)\leq\int_0^1 f\big(tx+(1-t)y\big)dt
                                +\varepsilon \qquad(x,y\in D).
}
and
\Eq{nrs}{
     \int_0^1 f\big(tx+(1-t)y\big)dt
                  \leq\frac{f(x)+f(y)}{2}+\varepsilon \qquad(x,y\in D).
}
Concerning the reversed implication, Nikodem, Riedel, and Sahoo
in \cite{NikRieSah07} have recently shown that the
$\varepsilon$-Hermite--Hadamard inequalities \eq{nrsl} and \eq{nrs}
do not imply the $c\varepsilon$-convexity of $f$ (with any $c>0$).
Thus, in order to obtain results that establish implications between the
approximate Hermite--Hadamard inequalities and the approximate Jensen inequality,
one has to consider these inequalities with nonconstant error terms.
More precisely, we will investigate the connection between the following
functional inequalities:
\Eq{J}{
  f\Big(\frac{x+y}{2}\Big)\leq\frac{f(x)+f(y)}{2}+\alpha_J(x-y) \qquad (x,y\in D),
}
\Eq{lH}{
  f\left(\frac{x+y}{2}\right)\leq \int_0^1 f\big(tx+(1-t)y\big)dt +
  \alpha_H(x-y) \qquad(x,y\in D),
}
\Eq{H}{
\int_0^1f\big(tx+(1-t)y\big)\rho(t)dt
\leq\lambda f(x)+(1-\lambda)f(y) +\alpha_H(x-y) \qquad (x,y\in D),
}
where $\alpha_H,\alpha_J:D^*\to\R_+$ are given even functions, $\lambda\in[0,1]$,
and $\rho:[0,1]\to\R_+$ is an integrable nonnegative function with
\Eq{*}{
  \int_0^1\rho(t)dt=1.
}

In order to describe the old and new results about the connection of the
approximate Jensen convexity inequality \eq{J} and the
approximate lower and upper Hermite--Hadamard inequalities \eq{lH} and \eq{H},
we need to introduce the following terminology.

For a function $f:D\to\R$, we say that $f$ is \textit{lower hemicontinuous},
\textit{upper hemicontinuous}, and \textit{hemiintegrable} on $D$
if, for all $x,y\in D$, the mapping
\Eq{*}
{
t\mapsto f(tx+(1-t)y) \qquad (t\in [0,1])
}
is lower semicontinuous, upper semicontinuous, and Lebesgue integrable on $[0,1]$, respectively.
We say that a function $h:D^*\to\R$ is \textit{radially lower semicontinuous},
\textit{radially upper semicontinuous}, \textit{radially increasing}, \textit{radially measurable},
and \textit{radially bounded}, if for all $u\in D^*$, the mapping
\Eq{hh}
{
t\mapsto h(tu) \qquad (t\in [0,1])
}
is lower semicontinuous, upper semicontinuous, increasing, measurable, and bounded on $[0,1]$,
respectively. (Note that, by the convexity of $D^*$ and $0\in D^*$, the function \eq{hh} is correctly defined.)

In \cite{HazPal09}, the relationships between the approximate lower
Hermite--Hadamard inequality \eq{lH} and approximate Jensen convexity inequality
\eq{J} were examined by H\'azy and P\'ales, who obtained the following results.

\THM{A}{Let $\alpha_J:D^*\to\R_+$ be a radially Lebesgue
integrable even function. Assume that $f:D\to\R$ is hemiintegrable
on $D$ and satisfies the approximate Jensen inequality \eq{J}.
Then $f$ also satisfies the approximate lower Hermite--Hadamard inequality
\eq{lH}, where $\alpha_H:D^*\to\R$ is defined by
\Eq{JH}{
  \alpha_H(u):=\int_0^1 \alpha_J(|1-2t|u)dt\qquad(u\in D^*).
}}

\THM{B}{Let $\alpha_H:D^*\to\R_+$ be an even function.
Assume that $f:D\to \R$ is an upper hemicontinuous function
satisfying the approximate lower Hermite--Hadamard inequality \eq{lH}.
Then $f$ satisfies the approximate Jensen inequality \eq{J} where
$\alpha_J:2D^*\to\R_+$ is a radially increasing nonnegative
solution of the functional inequality
\Eq{DI}{
  \int_0^{1}\alpha_J(2tu)dt+\alpha_H(u)\leq\alpha_J(u)\qquad (u\in D^*).
}}

The main new results of Section 2 are the following two theorems that are
analogous to \thm{A} above.

\Thm{1}{Let $\alpha_J:D^*\to\R$ be radially bounded, measurable and $\rho:[0,1]\to\R_+$ be
a Lebesgue integrable function with $\int_0^1\rho=1$.
Assume that $f:D\to\R$ is hemiintegrable and approximately Jensen convex in the sense of \eq{J}.
Then $f$ also satisfies the approximate upper Hermite--Hadamard inequality \eq{H}
with $\lambda:=\int_0^1 t\rho(t)dt$ and $\alpha_H:D^*\to\R$ defined by
\Eq{1h} {
\alpha_H(u):=\sum_{n=0}^\infty \frac{1}{2^n}
   \int_{0}^1\alpha_J\big(2d_{\Z}(2^{n}t)u\big)\rho(t)dt    \qquad (u\in D^*),
}
where, for $s\in\R$, $d_\Z(s):=\dist(s,\Z)=\inf\{|s-k|:k\in\Z\}$.}

\Thm{3}{Let $\alpha_J:D^*\to\R_+$ be radially increasing such that
\Eq{cc}{
\sum_{n=0}^{\infty}\alpha_J\Big(\frac{u}{2^n}\Big)<\infty\qquad(u\in D^*).
}
If $f : D \to \R$ is upper hemicontinuous and $\alpha_J$-Jensen
convex on $D$, i.e., \eq{J} holds, then $f$ also satisfies the Hermite--Hadamard inequality \eq{H}
with $\lambda:=\int_0^1 t\rho(t)dt$ and $\alpha_H:D^*\to\R$ is defined by
\Eq{TTb}
{
\alpha_H(u):=\sum_{n=0}^{\infty}2\alpha_J\Big(\frac{u}{2^n}\Big)\int_0^1d_{\Z}(2^nt)\rho(t)dt
     \qquad (u\in D^*).
}}

The main result of Section 3 is the following theorem which
corresponds to \thm{B} above.

\Thm{2}{Let $\alpha_H:D^*\to\R$ be even and radially upper semicontinuous,
$\rho:[0,1]\to\R_+$ be integrable with $\int_0^1\rho=1$ and there exist $c\geq0$ and
$p>0$ such that
\Eq{c}{
  \rho(t)\leq c(-\ln|1-2t|)^{p-1} \qquad(t\in ]0,\tfrac12[\cup]\tfrac12,1[),
}
and $\lambda\in[0,1]$. Then every $f:D\to\R$ lower hemicontinuous function
satisfying the approximate upper Hermite--Hadamard inequality \eq{H},
fulfills the approximate Jensen inequality \eq{J} provided that
$\alpha_J:D^*\to\R$ is a radially lower semicontinuous solution of
the functional inequality
\Eq{a} {
\alpha_J(u)\geq
\int_0^1\alpha_J(|1-2t|u)\rho(t)dt+\alpha_H(u)\qquad (u\in D^*)
}
and $\alpha_J(0)\geq\alpha_H(0)$.}

In Section 2, implications from inequality \eq{J} to \eq{H} will be investigated.

A weaker form of \thm{1} could be deduced from the following result which 
was obtained by the authors in \cite{MakPal10b}.

\THM{MP}{Let $\alpha_J:D^*\to\R_+$ be radially bounded and even.
Then, an upper hemicontinuous function $f:D\to\R$ is $\alpha_J$-Jensen convex on $D$,
i.e., \eq{J} holds if and only if
\Eq{TTd} {
  f(tx+(1-t)y)\leq tf(x)+(1-t)f(y)+\sum_{n=0}^{\infty}\frac{\alpha_J\big(2d_{\Z}(2^nt)(x-y)\big)}{2^n}
  \qquad(x,y\in D,\,t\in [0,1]).
}}

In the proof of \thm{1}, we will directly derive \eq{H} from \eq{J} under more general
circumstances.

To deduce \thm{3}, the following result of Jacek and J\'ozef Tabor \cite{TabTab09b}
will be used.

\THM{TT}{
Let $\alpha_J:D^*\to\R_+$ be radially increasing and even such that, for all $u\in D^*$, 
$\sum_{n=0}^{\infty}\alpha_J\big(2^{-n}u\big)<\infty$.
Then, an upper hemicontinuous function $f:D\to\R$ is $\alpha_J$-Jensen convex on $D$,
i.e., \eq{J} holds if and only if
\Eq{TTa} {
  f(tx+(1-t)y)\leq tf(x)+(1-t)f(y)+\sum_{n=0}^{\infty}2\alpha_J\Big(\frac{x-y}{2^n}\Big)d_{\Z}(2^nt)
  \qquad(x,y\in D,\,t\in [0,1]).
}}

In the particular case when $\alpha_J$ is a linear combination of power
functions, we also deduce some consequences of \thm{1} and \thm{3}.
For this aim, we will have to recall two different notions of Takagi type functions.
For $q>0$, define the functions $T_q:\R\to\R$ and $S_q:\R\to\R$ by
\Eq{hp}{
  T_q(x):=\sum_{n=0}^\infty\frac{\big(2d_\Z(2^nx)\big)^q}{2^n},\qquad
  S_q(x):=\sum_{n=0}^\infty\frac{d_\Z(2^nx)}{2^{nq-1}}\qquad(x\in\R).
}
They generalize the classical Takagi function $T_1=S_1=T$ in two ways. 
It is more difficult to see that $T_2=S_2$ is also valid.
These functions have an important role in approximate convex analysis.

The importance of the functions $T_q$ introduced above is
enlightened by the following result (cf.\ \cite{HazPal04}, \cite{HazPal05},
\cite{Haz05a}, \cite{Haz07a}) which is a generalization of
the celebrated Bernstein--Doetsch theorem \cite{BerDoe15}.

\THM{H7}{Let $X$ be a normed space, $q>0$ and $a\geq0$.
Then a locally upper bounded function $f:D \to \R$ is $(a,q)$-Jensen
convex on $D$, i.e.,
\Eq{*}{
  f\Big(\frac{x+y}{2}\Big)\leq\frac{f(x)+f(y)}{2}+a\|x-y\|^q \qquad (x,y\in D),
}
if and only
\Eq{H7a} {
  f(tx+(1-t)y)\leq tf(x)+(1-t)f(y)+a T_q(t)\|x-y\|^q \qquad  (x,y\in D,\,t\in[0,1]).
} }

The other Takagi type function $S_q$ was introduced by
Tabor and Tabor. Its role and importance
in the theory of approximate convexity is shown by the next theorem
(\cite{TabTab09b}, \cite{TabTab09a}).

\THM{TTa}{Let $X$ be a normed space, $q>0$ and $a\geq0$.
Then a locally upper bounded function $f:D \to \R$ is $(a,q)$-Jensen
convex on $D$, i.e.,
\Eq{*}{
  f\Big(\frac{x+y}{2}\Big)\leq\frac{f(x)+f(y)}{2}+a\|x-y\|^q \qquad (x,y\in D),
}
if and only if
\Eq{TTc} {
  f(tx+(1-t)y)\leq tf(x)+(1-t)f(y)+aS_q(t)\|x-y\|^q \qquad  (x,y\in D,\,t\in[0,1]).
}}
In view of the results in the papers \cite{MakPal10b} and \cite{MakPal12b}, the error terms 
in \eq{H7a} and \eq{TTc} are the best possible if $0<q\leq1$ and $1\leq q\leq 2$, respectively.

In Section 3, for every parameter $p>0$, we define a class of functions, denoted by $\Phi_p$. A certain 
convolution-like operation is also introduced in $\bigcup_{p>0}\Phi_p$ and its properties are 
described in Propositions \ref{Ppp}--\ref{Pgphi}. 
These tools will be instrumental in the proof of \thm{2} which will be carried out in several steps.
Finally, when the error function $\alpha_H$ is a linear combination of power functions,
we will also deduce some corollaries of \thm{2}.

\section{From Jensen inequality to Hermite--Hadamard inequality}
\setcounter{theorem}{0}

The following statement will be essential to obtain our first main result, \thm{1}.

\Prp{ppp}{
Let $\rho:[0,1]\to\R_+$ be a Lebesgue integrable function.
Then, the function $\psi:[0,1]\to\R$ defined by
\Eq{Psi}
{
\psi(t)=\frac12\sum_{n=0}^{\infty}\frac{1}{4^n}
   \Big(\sum_{k=0}^{2^n-1}\rho(\tfrac{t+k}{2^n})\Big)\qquad (t\in[0,1]),
}
is a nonnegative integrable solution of the functional equation
\Eq{phi} {
\rho(t)=2\psi(t)-\frac{\psi(\tfrac{t}{2})+\psi(\tfrac{t+1}{2})}2\qquad
(t\in[0,1]).
}
Furthermore,
\Eq{lambda}
{
\int_0^1\psi=\int_0^1\rho,\qquad
\int_{\frac12}^1\psi=\int_0^1 t\rho(t)dt \qquad \mbox{and} \qquad
\int_0^{\frac12}\psi= \int_0^1 (1-t)\rho(t)dt.
}}

\begin{proof}
Define the sequence $\psi_n:[0,1]\to\R$, by
\Eq{psi}{
\psi_0:=\tfrac12\rho, \qquad
\psi_n(t):=\frac12\rho(t)+\frac14\big(\psi_{n-1}(\tfrac{t}{2})+\psi_{n-1}(\tfrac{t+1}{2})\big)
\qquad (t\in[0,1],\,n\in\N).
}
Then the sequence $(\psi_n)$ is nondecreasing, i.e.,
\Eq{n}
{
0\leq\psi_{n-1}\leq \psi_{n},\qquad \mbox{and} \qquad \int_0^1\psi_n=\frac{2^{n+1}-1}{2^{n+1}}\int_0^1\rho \qquad  (n\in\N).
}
We prove \eq{n} by induction on $n\in\N$.
For $n=1$, by the definition of $\psi_1$ and the nonnegativity of $\rho$, for $t\in[0,1]$, we have that
\Eq{*}
{
\psi_{1}(t)=\frac12\rho(t)+\frac14\big(\psi_0(\tfrac{t}{2})+\psi_0(\tfrac{t+1}{2})\big)
           =\frac12\rho(t)+\frac18\big(\rho(\tfrac{t}{2})+\rho(\tfrac{t+1}{2})\big)
           \geq \frac12\rho(t)=\psi_0(t),
}
and
\Eq{*}
{
\int_0^1\psi_0=\frac{1}{2}\int_0^1\rho.
}
Assume that, for some $n\in\N$, \eq{n} holds and consider the case $n+1$.
By the definition of $\psi_{n+1}$, the inductive assumption and the nonnegativity of $\psi_n$,
for $t\in[0,1]$, yields
\Eq{*}
{
\psi_{n+1}(t)=\frac12\rho(t)+\frac14\big(\psi_n(\tfrac{t}{2})+\psi_n(\tfrac{t+1}{2})\big)
           \geq\frac12\rho(t)+\frac14\big(\psi_{n-1}(\tfrac{t}{2})+\psi_{n-1}(\tfrac{t+1}{2})\big)
           =\psi_n(t).
}
Using the definition $\psi_{n+1}$, the substitution $s:=\frac{t}{2}$ and $s:=\frac{t+1}{2}$,
finally the inductive assumption, we get
\Eq{*}{
\int_0^1\psi_{n+1}
       &=\frac12\int_0^1\rho+\frac14\Big(\int_0^1\psi_n(\tfrac{t}{2})dt+\int_0^1\psi_n(\tfrac{t+1}{2})dt\Big)
       =\frac12\int_0^1\rho+\frac12\Big(\int_0^{\tfrac12}\psi_n+\int_{\tfrac12}^1\psi_n\Big)\\
       &=\frac12\int_0^1\rho+\frac12\int_0^{1}\psi_n
       =\Big(\frac12+\frac{2^{n+1}-1}{2^{n+2}}\Big)\int_0^1\rho=\frac{2^{n+2}-1}{2^{n+2}}\int_0^1\rho.
}
Denote by $L^1[0,1]$ the space of Lebesgue integrable functions $\chi:[0,1]\to\R$.
Then $L^1[0,1]$ is a Banach-space with the standard norm $\|\chi\|_1:=\int_0^1|\chi|$.
Now we prove that $(\psi_n)$ is a Cauchy sequence in $L^1[0,1]$. Using \eq{n},
for $n\leq m$, we get that
\Eq{*}{
 \|\psi_m-\psi_n\|_1=
 \int_0^1(\psi_m-\psi_n)=\frac{2^m-2^n}{2^{n+m+1}}\int_0^1\rho
             \leq\frac{1}{2^n}\int_0^1\rho,
}
which implies that $(\psi_n)$ is indeed a Cauchy sequence. Hence it converges to
a function $\psi\in L^1[0,1]$.
To prove \eq{Psi}, we show, by induction on $n\in\N$, that
\Eq{Psi_n}
{
\psi_n(t)=\frac12\sum_{i=0}^{n}\frac{1}{4^i}
   \Big(\sum_{k=0}^{2^i-1}\rho(\tfrac{t+k}{2^i})\Big)\qquad (t\in[0,1])
}
holds. For $n=0$, we have an obvious identity.
Assume that \eq{Psi_n} holds some $n\in\N$.
Using the definition of $\psi_{n+1}$ and the inductive assumption, we obtain
\Eq{*}
{
\psi_{n+1}(t)=&\frac12\rho(t)+\frac14\big(\psi_{n}(\tfrac{t}{2})+\psi_{n}(\tfrac{t+1}{2})\big)\\
             =&\frac12\rho(t)+\frac18\sum_{i=0}^{n}\frac{1}{4^i}
                       \Big(\sum_{k=0}^{2^i-1}\rho(\tfrac{t+2k}{2^{i+1}})\Big)
                             +\frac18\sum_{i=0}^{n}\frac{1}{4^i}
                       \Big(\sum_{k=0}^{2^i-1}\rho(\tfrac{t+2k+1}{2^{i+1}})\Big)\\
             =&\frac12\rho(t)+\frac12\sum_{i=0}^{n}\frac{1}{4^{i+1}}
                       \sum_{k=0}^{2^i-1}\Big(\rho(\tfrac{t+2k}{2^{i+1}})
                             +\rho(\tfrac{t+2k+1}{2^{i+1}})\Big)\\
             =&\frac12\rho(t)+\frac12\sum_{i=0}^{n}\frac{1}{4^{i+1}}
                       \sum_{k=0}^{2^{i+1}-1}\rho(\tfrac{t+k}{2^{i+1}})
             =\frac12\sum_{i=0}^{n+1}\frac{1}{4^{i}}
                       \Big(\sum_{k=0}^{2^{i}-1}\rho(\tfrac{t+k}{2^{i}})\Big),
}
which proves \eq{Psi_n}.
Thus, taking the limit $n\to\infty$ in \eq{Psi_n}, we obtain \eq{Psi}.

To prove the first expression in \eq{lambda}, integrate \eq{phi} on $[0,1]$, then we get
\Eq{*}
{
\int_0^1\rho&=2\int_0^1\psi-\frac{\int_0^1\psi(\tfrac{t}{2})dt+\int_0^1\psi(\tfrac{t+1}{2})dt}2
                 =2\int_0^1\psi-\Big(\int_0^{\frac12}\psi+\int_{\frac12}^1\psi\Big)
                  =\int_0^1\psi.
}
To prove the second expression in \eq{lambda}, multiply \eq{phi} by $t$
and integrate it on $[0,1]$. Thus we get
\Eq{*}
{
\int_0^1t\rho(t)dt
&=2\int_0^1t\psi(t)dt-\Big(\int_0^1\tfrac{t}{2}\psi(\tfrac{t}{2})dt
                       +\int_0^1\tfrac{t+1}{2}\psi(\tfrac{t+1}{2})dt\Big)+\frac12\int_0^1\psi(\tfrac{t+1}{2})dt\\
&=2\int_0^1t\psi(t)dt-2\Big(\int_0^{\frac12}s\psi(s)ds
                       +\int_{\frac12}^1s\psi(s)ds\Big)+\int_{\frac12}^1\psi=\int_{\frac12}^1\psi.
}
The last equality in \eq{lambda} is a consequence of the first and second equalities.
\end{proof}

\begin{proof}[Proof of \thm{1}]
By \prp{ppp}, the function $\psi:[0,1]\to\R_+$ defined by \eq{Psi} is a 
Lebesgue integrable function satisfying the functional equation \eq{phi} for which
\eq{lambda} holds.

Let $f:D\to\R$ be an approximately Jensen convex function. Let $x,y\in D$ be arbitrary fixed.
Then, by approximate Jensen convexity of $f$, we have that
\Eq{*} {
f(tx+(1-t)y)\leq
\begin{cases}
\dfrac{f(2tx+(1-2t)y)+f(y)}{2}+\alpha_J\big(2t(x-y)\big)\qquad &(t\in[0,\tfrac12]),\\[2mm]
\dfrac{f(x)+f((2t-1)x+(2-2t)y)}{2}+\alpha_J\big((2-2t)(x-y)\big)\qquad &(t\in[\tfrac12,1]).
\end{cases}
}
Multiplying the above inequality by $2\psi(t)$,
taking the integral over $[0,1]$, we get
\Eq{1a} {
\int_{0}^{1}f(tx&+(1-t)y)2\psi(t)dt
            \leq\int_{0}^{\frac12}\hspace{-2mm}\Big(f(2tx+(1-2t)y)+f(y)
                   +2\alpha_J\big(2t(x-y)\big)\Big)\psi(t)dt\\
             &+\int_{\frac12}^1\hspace{-2mm}\Big(f(x)+f((2t-1)x+(2-2t)y)
                  +2\alpha_J\big((2-2t)(x-y)\big)\Big)\psi(t)dt.
}
Substituting $t:=\frac{s}{2}$ and $t:=\frac{1+s}{2}$ in the first and second terms
on the right hand side of \eq{1a}, using \eq{lambda}, and observing that
$2d_\Z(t)=\min(2t,2-2t)$, we have that
\Eq{1c}
{
\int_{0}^{\frac12}&\hspace{-1mm}\Big(f(2tx+(1-2t)y)+f(y)
                   +2\alpha_J\big(2t(x-y)\big)\Big)\psi(t)dt\\
            &=f(y)\int_{0}^1 (1-t)\rho(t)dt
              +\frac12\int_{0}^1f(sx+(1-s)y)\psi(\tfrac{s}{2})ds
              +2\int_{0}^{\frac12}\alpha_J\big(2d_{\Z}(t)(x-y)\big)\psi(t)dt,\\
\int_{\frac12}^1&\hspace{-1mm}\Big(f(x)+f((2t-1)x+(2-2t)y)
                  +2\alpha_J\big((2-2t)(x-y)\big)\Big)\psi(t)dt\\
       &=f(x)\int_0^1 t\rho(t)dt
         +\frac12\int_{0}^1f(sx+(1-s)y)\psi(\tfrac{1+s}{2})ds
         +2\int_{\frac12}^1\alpha_J\big(2d_{\Z}(t)(x-y)\big)\psi(t)dt.
}
Combining \eq{1a} and \eq{1c}, we get that
\Eq{*}
{
\int_{0}^{1}f(tx&+(1-t)y)\rho(t)dt
=\int_{0}^{1}f(tx+(1-t)y)\Big(2\psi(t)-\tfrac12\psi(\tfrac{t}2)-\tfrac12\psi(\tfrac{t+1}2)\Big)dt\\
            &\leq f(x)\int_{0}^1 t\rho(t)dt+f(y)\int_{0}^1 (1-t)\rho(t)dt
             +2\int_{0}^1\alpha_J\big(2d_{\Z}(t)(x-y)\big)\psi(t)dt.
}
To complete the proof, it remains to show that the last term containing $\psi$
equals $\alpha_H(x-y)$. Indeed, applying formula \eq{Psi} and the 1-periodicity
of the function $d_\Z$, we get
\Eq{*}{
  2\int_{0}^1\alpha_J\big(2d_{\Z}(t)(x-y)\big)\psi(t)dt
  &=2\int_{0}^1\alpha_J\big(2d_{\Z}(t)(x-y)\big)\frac12\sum_{n=0}^{\infty}\frac{1}{4^n}
   \Big(\sum_{k=0}^{2^n-1}\rho(\tfrac{t+k}{2^n})dt\Big)\\
  &=\sum_{n=0}^{\infty}\frac{1}{4^n}
   \Big(\sum_{k=0}^{2^n-1}\int_{0}^1\alpha_J\big(2d_{\Z}(t)(x-y)\big)\rho(\tfrac{t+k}{2^n})dt\Big)\\
  &=\sum_{n=0}^{\infty}\frac{1}{4^n}
   \Big(\sum_{k=0}^{2^n-1}2^n\int_{\frac{k}{2^n}}^{\frac{k+1}{2^n}}
              \alpha_J\big(2d_{\Z}(2^ns-k)(x-y)\big)\rho(s)ds\Big)\\
  &=\sum_{n=0}^{\infty}\frac{1}{2^n}
   \int_0^1\alpha_J\big(2d_{\Z}(2^ns)(x-y)\big)\rho(s)ds=\alpha_H(x-y),
}
which proves the statement.
\end{proof}

The other form of the error function $\alpha_H$ stated in \thm{3} can be obtained by using
\thm{TT} by Jacek Tabor and J\'ozef Tabor \cite{TabTab09b}.

\begin{proof}[Proof of \thm{3}]
If $f : D \to \R$ is upper hemicontinuous and $\alpha_J$-Jensen
convex on $D$, then \eq{TTa} holds. Multiplying this inequality by $\rho(t)$ and then
integrating with respect to $t$ over $[0,1]$, \eq{TTb} follows immediately.
\end{proof}

Let $X$ be a normed space. Next, we consider the case, when
$\alpha_J$ is a linear combination of the powers
of the norm with positive exponents, i.e., if $\alpha_J$ is of the form
\Eq{DJP+}{
  \alpha_J(u):=\int\limits_{]0,\infty[} \|u\|^qd\mu_J(q)\qquad (u\in D^*),
}
where $\mu_J$ is a nonnegative Borel measure on the interval $]0,\infty[$.
An important particular case is when $\mu_J$ is of the form $\sum_{i=1}^k c_i\delta_{q_i}$,
where $c_i\in\R_+$, $q_i>0$ and $\delta_{q_i}$ stands for the Dirac measure concentrated at $q_i$ 
for $i\in\{1,\dots,k\}$.

\Thm{A2}{Let $\rho:[0,1]\to\R_+$ be a Lebesgue integrable function with $\int_0^1\rho=1$
and let $\mu_J$ be a signed Borel measure on $]0,\infty[$ such that
\Eq{*}
{
  \int\limits_{]0,\infty[}\|u\|^qd|\mu_J|(q)<\infty\qquad (u\in D^*).
}
Assume that $f:D\to\R$ is hemiintegrable on $D$
and is approximately Jensen convex in the following sense
\Eq{A2a}
{
f\Big(\frac{x+y}{2}\Big)\leq \frac{f(x)+f(y)}{2}
  +\int\limits_{]0,\infty[}\|x-y\|^qd\mu_J(q)\qquad (x,y\in D).
}
Then $f$ also satisfies the approximate Hermite--Hadamard inequality
\Eq{A2b}
{
\int_{0}^{1}\!\!f(tx+(1-t)y)\rho(t)dt\leq \lambda f(x)+(1-\lambda)f(y)
  +\!\!\int\limits_{]0,\infty[}\!\!\int_0^1\!\!T_q(t)\rho(t)dt\|x-y\|^q d\mu_J(q)\quad (x,y\in D),
}
with $\lambda:=\int_0^1 t\rho(t)dt$.}

\begin{proof}
It is easy to see that $\alpha_J$ defined by \eq{DJP+} is radially bounded and measurable.
Thus, by \thm{1}, it is enough to compute the error function $\alpha_{H}$ defined by \eq{1h}.
Hence, using  \eq{1h}, \eq{DJP+}, Fubini's theorem and Lebesgue's theorem, we obtain
\Eq{*}
{
\alpha_H(u)&=\sum_{n=0}^\infty \frac{1}{2^n}\int_{0}^1
             \!\!\int\limits_{]0,\infty[} \|\big(2d_{\Z}(2^{n}t)u\big)\|^qd\mu_J(q)\rho(t)dt\\
           &= \int\limits_{]0,\infty[}\!\!\int_{0}^1 \sum_{n=0}^\infty
                 \frac{\big(2d_{\Z}(2^{n}t)\big)^q}{2^n}\rho(t)dt\|u\|^qd\mu_J(q)
           =\int\limits_{]0,\infty[}\!\!\int_{0}^1 T_q(t)\rho(t)dt\|u\|^qd\mu_J(q) \quad (u\in D^*),
}
which completes the proof.
\end{proof}

\Thm{A2+}{Let $\rho:[0,1]\to\R_+$ be a Lebesgue integrable function with $\int_0^1\rho=1$
and let $\mu_J$ be a nonnegative Borel measure on $]0,\infty[$, such that
\Eq{cc1}{
\int\limits_{]0,\infty[}\|u\|^qd\mu_J(q)<\infty \qquad (u\in D^*)
}
and
\Eq{cc2}{
\int\limits_{]0,\infty[}\frac{2^q}{2^q-1}d\mu_J(q)<\infty.
}
Assume that $f:D\to\R$ is upper hemicontinuous and approximately Jensen convex in the sense of \eq{A2a}.
Then $f$ also satisfies the following approximate Hermite--Hadamard inequality
\Eq{Sq}
{
\int_{0}^{1}\!\!\!f(tx+(1-t)y)\rho(t)dt\!\leq \lambda f(x)+(1-\lambda)f(y)
  +\!\!\!\int\limits_{]0,\infty[}\!\!\!\int_0^1\!\!S_q(t)\rho(t)dt\|x-y\|^q d\mu_J(q)\quad (x,y\in D),
}
with $\lambda:=\int_0^1 t\rho(t)dt$.}

\begin{proof}
Consider the function $\alpha_J$ defined by \eq{DJP+}.
Then, for all $u\in D^*$, the mapping $t\mapsto\alpha_J(tu)=\int\limits_{]0,\infty[}t^q\|u\|^qd\mu_J(q)$
is increasing on $[0,1]$ and, for all $u\in D^*$,
\Eq{*}{
\sum_{n=0}^{\infty}\alpha_J\Big(\frac{u}{2^n}\Big)
       =\sum_{n=0}^{\infty}\int\limits_{]0,\infty[}\|2^{-n}u\|^qd\mu_J(q)
       =\int\limits_{]0,\infty[}\frac{2^q}{2^q-1}\|u\|^qd\mu_J(q).
}
If $\|u\|\leq 1$, then the latter series is convergent in virtue of \eq{cc2}.
For $\|u\|>1$, we have
\Eq{*}{
\int\limits_{]0,\infty[}\frac{2^q}{2^q-1}\|u\|^qd\mu_J(q)
  \leq \|u\|\int\limits_{]0,1]}\frac{2^q}{2^q-1}d\mu_J(q)
    +2\int\limits_{]1,\infty[}\|u\|^qd\mu_J(q)<\infty,
}
which proves the convergence condition \eq{cc}.
Thus, by \thm{3}, it is enough to compute the error function $\alpha_{H}$ defined by \eq{TTb}.
Hence, using  \eq{TTb}, \eq{DJP+}, Fubini's theorem and Lebesgue's theorem, we obtain
\Eq{*}
{
\alpha_H(u)&=\sum_{n=0}^{\infty}2\int\limits_{]0,\infty[}\Big(\frac{\|u\|}{2^n}\Big)^qd\mu_J(q)\int_0^1d_{\Z}(2^nt)\rho(t)dt\\
           &=\int\limits_{]0,\infty[}\!\!\int_{0}^1 \sum_{n=0}^\infty \frac{d_{\Z}(2^nt)}{2^{nq-1}}\rho(t)dt\|u\|^qd\mu_J(q)
           =\int\limits_{]0,\infty[}\!\!\int_{0}^1 S_q(t)\rho(t)dt\|u\|^qd\mu_J(q) \quad (u\in D^*),
}
which completes the proof.
\end{proof}

Now we consider the case in the previous theorems when
$\rho\equiv1$ and the measure $\mu_J$ is the Dirac
measure $a\delta_q$.

\Cor{A3}{Let $a\in\R_+$ and $q>0$. Assume that $f:D\to\R$ is hemiintegrable
and satisfies the following approximate Jensen convexity inequality
\Eq{A3a}
{
f\Big(\frac{x+y}{2}\Big)\leq \frac{f(x)+f(y)}{2}+a\|x-y\|^q\qquad (x,y\in D).
}
Then $f$ also satisfies the following approximate Hermite--Hadamard inequality
\Eq{A3b}
{
\int_{0}^{1}f(tx+(1-t)y)dt\leq \frac{f(x)+f(y)}{2}+\frac{2a}{q+1}\|x-y\|^q\qquad (x,y\in D).
}
}
\begin{proof}
The conditions of \thm{A2} hold with $\rho\equiv1$ and $\mu_J:=a\delta_q$,
Then \eq{A2a} holds by \eq{A3a}, hence
to prove the statement, it is enough to compute the error term in \eq{A2b}.
Using the definition of the $T_q$, the substitution $s:=2^nt$ and the
$1$-periodicity of $d_\Z^q$, we get
\Eq{*}
{
\int_0^1T_q(t)dt&=\int_0^1\bigg(\sum_{n=0}^\infty\frac{\big(2d_\Z(2^nt)\big)^q}{2^n}\bigg)dt
                 =\sum_{n=0}^\infty\frac{2^q}{2^n}\int_0^1\big(d_\Z(2^nt)\big)^qdt
                 =\sum_{n=0}^\infty\frac{2^q}{2^{2n}}\int_0^{2^n}\big(d_\Z(s)\big)^qds\\
                &=\sum_{n=0}^\infty\frac{2^q}{2^{n}}\int_0^{1}\big(d_\Z(s)\big)^qds
                 =\sum_{n=0}^\infty\frac{2^q}{2^{n}}\bigg(\int_0^{\frac12}s^qds+\int_{\frac12}^1(1-s)^qds\bigg)\\
                 &=\sum_{n=0}^\infty\frac{2^q}{2^{n}}\Big(\frac{(\tfrac12)^{q+1}}{q+1}+\frac{(\tfrac12)^{q+1}}{q+1}\Big)
                 =\frac1{q+1}\sum_{n=0}^\infty\frac1{2^n}=\frac2{q+1} \qquad (u\in D^*).
                 }
Thus, \eq{A2b} reduces to \eq{A3b}, which proves the statement.
\end{proof}

\Cor{A3+}{Let $a\in\R_+$ and $q>0$.
Assume that $f:D\to\R$ is upper hemicontinuous and satisfies the approximate Jensen convexity inequality \eq{A3a}.
Then $f$ also satisfies the following approximate Hermite--Hadamard inequality
\Eq{A3b+}
{
\int_{0}^{1}f(tx+(1-t)y)dt\leq \frac{f(x)+f(y)}{2}+\frac{2^qa}{2^{q+1}-2}\|x-y\|^q   \qquad (x,y\in D).
}}

\begin{proof}
The conditions of \thm{A2+} are satisfied with $\rho\equiv1$ and $\mu_J:=a\delta_q$,
Then \eq{A2a} holds by \eq{A3a}, hence
to prove the statement, it is enough to compute the error term in \eq{Sq}.
Using the definition of the $S_q$, the substitution $s:=2^nt$ and the $1$-periodicity of $d_\Z$, we get
\Eq{*}
{
\int_0^1S_q(t)dt&=\int_0^1\bigg(\sum_{n=0}^\infty\frac{d_\Z(2^nt)}{2^{nq-1}}\bigg)dt
                 =\sum_{n=0}^\infty\frac{1}{2^{nq-1}}\int_0^1d_\Z(2^nt)dt
                 =\sum_{n=0}^\infty\frac{1}{2^{nq+n-1}}\int_0^{2^n}d_\Z(s)ds\\
                &=\sum_{n=0}^\infty\frac{1}{2^{nq-1}}\int_0^{1}d_\Z(s)ds
                 =\sum_{n=0}^\infty\frac{1}{2^{nq-1}}\bigg(\int_0^{\frac12}sds+\int_{\frac12}^1(1-s)ds\bigg)\\
                &=\sum_{n=0}^\infty\frac{1}{2^{nq-1}}\frac{1}{4}
                 =\frac12\sum_{n=0}^\infty\frac{1}{2^{nq}}=\frac{2^q}{2^{q+1}-2} \qquad (u\in D^*).
                 }
Thus, \eq{Sq} reduces to \eq{A3b+}, which completes the proof.
\end{proof}

\Rem{*}{The constants in the two error terms obtained in \eq{A3b} and \eq{A3b+}
are comparable in the following way: for $q\in]0,1[\cup]2,\infty[$,
\Eq{*}{
  \frac{2}{q+1}<\frac{2^q}{2^{q+1}-2}
}
and the inequality reverses for $q\in]1,2[$.}

\section{From Hermite--Hadamard inequality to Jensen inequality}
\setcounter{theorem}{0}

For $p>0$, define the class of functions $\Phi_p$ by
\Eq{*}
{\Phi_p:=\Big\{\phi:]0,1[\to \R\mid \phi \mbox{ is Lebesgue
measurable and }
   \/\phi\/_{p}:=\sup_{t\in]0,1[}|\ln t|^{1-p}|\phi(t)|<\infty\Big\}.
}

In the sequel, $\Gamma$ denotes Euler's Gamma function.

\Prp{pp}{For all $p>0$, the elements of $\Phi_p$ are Lebesgue
integrable functions and
\Eq{L1}{
   \|\phi\|_1=\int_0^1|\phi(t)|dt\leq \Gamma(p)\/\phi\/_p\qquad (\phi\in\Phi_p).
}}

\begin{proof}
Let $p>0$ and $\phi\in \Phi_p$. From the definition of $\Phi_p$, we get that
 \Eq{*}{
 |\phi(t)|\leq \/\phi\/_p(-\ln t)^{p-1}\qquad (t\in]0,1[).
 }
 Thus, with the substitution $s=-\ln t$, we get
\Eq{*}{
  \int_0^1|\phi(t)|dt\leq \/\phi\/_p\int_0^1(-\ln t)^{p-1}dt
  =\/\phi\/_p \int_0^\infty s^{p-1}e^{-s}ds=\Gamma(p)\/\phi\/_p<\infty,
}
which proves the integrability of $\phi$ and \eq{L1}.
\end{proof}

\Prp{pq}{For $p,q>0$ and $\phi\in\Phi_p$, $\psi\in\Phi_q$, the
function $\phi*\psi$ defined by
\Eq{*} {
(\phi*\psi)(t):=\int_{t}^1\tfrac{1}{\tau}\phi(\tfrac{t}{\tau})\psi(\tau)d\tau\qquad
           (t\in]0,1[).
           }
is continuous on the open interval $]0,1[$, belongs to $\Phi_{p+q}$ and
\Eq{p+q}{
 \/\phi*\psi\/_{p+q}\leq
\frac{\Gamma(p)\Gamma(q)}{\Gamma(p+q)} \/\phi\/_{p}\/\psi\/_{q}. }
Furthermore,
\Eq{prod}{
\int_0^1(\phi*\psi)=\int_0^1\phi\int_0^1\psi. }}

\begin{proof}
Given $p,q>0$, it is well known that the function
\Eq{beta}{
   \tau\mapsto (1-\tau)^{p-1}\tau^{q-1}, \qquad(\tau\in]0,1[)
}
is integrable over $[0,1]$ and
\Eq{*}{
   B(p,q):=\int_0^1(1-\tau)^{p-1}\tau^{q-1}d\tau
     =\frac{\Gamma(p)\Gamma(q)}{\Gamma(p+q)}.
}

By the inclusions $\phi\in\Phi_p$, $\psi\in\Phi_q$, we have that
\Eq{pq}{
|\phi(t)|\leq \/\phi\/_p(-\ln t)^{p-1}\qquad\mbox{and}\qquad
|\psi(t)|\leq \/\psi\/_q(-\ln t)^{q-1}\qquad (t\in]0,1[). }

To prove that $\phi*\psi$ is continuous at $t\in]0,1[$, let $\eps>0$.
By the integrability of \eq{beta}, there exists $\rho\in]0,1[$ such
that, for every measurable subset $T\subseteq[0,1]$ with
$\meas(T)<\rho$,
\Eq{T}{
  \int_T(1-\tau)^{p-1}\tau^{q-1}d\tau<\frac{\eps}{6\/\phi\/_p\/\psi\/_q(-\ln t)^{p+q-1}+1}.
}
Define the function $\widetilde{\phi}:]0,\infty[\to\R$ by
\Eq{*} {
\widetilde{\phi}(x):=\phi(e^{-x}). 
}
Then, by the first inequality in \eq{pq}, we have that
\Eq{*} { x^{1-p}|\widetilde{\phi}(x)|\leq
\/\phi\/_p \qquad(x\in]0,\infty[). 
}
Applying Luzin's theorem for the bounded measurable function $x\mapsto
x^{1-p}\widetilde{\phi}(x)$, we can construct a measurable set
$\widetilde{H}\subseteq]0,\infty[$ and a continuous function
$\widetilde{f}:]0,\infty[\to\R$ such that
\Eq{w}{
  \meas(]0,\infty[\setminus \widetilde{H})  < \frac{(-\ln t)\rho}{2},\quad
  \widetilde{f}|_{\widetilde{H}}=\widetilde{\phi}|_{\widetilde{H}}
    \quad \mbox{and}\quad |\widetilde{f}(x)|\leq \/\phi\/_px^{p-1} \quad(x\in]0,\infty[).
}
Now define $f:]0,1[\to \R$ and $H\subseteq]0,1[$ by
\Eq{*} {
f(t):=\widetilde{f}(-\ln t)\qquad (t\in]0,1[)\qquad
\mbox{and}\qquad H:=\exp(-\widetilde{H}). }
In a view of \eq{w}, we get that
\Eq{fp}{
  f|_H=\phi|_H  \qquad \mbox{and}\qquad |f(t)|\leq \/\phi\/_p(-\ln t)^{p-1} \qquad(t\in]0,1[).
}
By the continuity of the logarithmic function, there exists
$\delta\in]0,\min(t,1-t)[$, such that, for all $s\in]t-\delta,t+\delta[$,
\Eq{ln}{
\/\phi\/_p\/\psi\/_q(-\ln(s))^{p+q-1}&<\/\phi\/_p\/\psi\/_q(-\ln(t))^{p+q-1}+\frac16
\qquad \mbox{and}\qquad\Big|1-\frac{\ln t}{\ln s}\Big|&<\rho.
}
For $s\in]0,1[$, we have
\Eq{II+}{
  |(\phi*\psi)(t)&-(\phi*\psi)(s)|
    =\Big|\int_{t}^1\tfrac{1}{\tau}\phi(\tfrac{t}{\tau})\psi(\tau)d\tau
      -\int_{s}^1\tfrac{1}{\tau}\phi(\tfrac{s}{\tau})\psi(\tau)d\tau\Big|\\
    &\leq\begin{cases}\displaystyle
     \int_s^t\tfrac{1}{\tau}\big|\phi(\tfrac{s}{\tau})\big||\psi(\tau)|d\tau
      +\int_{t}^1\tfrac{1}{\tau}\big|\phi(\tfrac{t}{\tau})
                   -\phi(\tfrac{s}{\tau})\big||\psi(\tau)|d\tau &\mbox{if } s<t,\\[3mm]\displaystyle
     \int_t^s\tfrac{1}{\tau}\big|\phi(\tfrac{t}{\tau})\big||\psi(\tau)|d\tau
      +\int_s^1\tfrac{1}{\tau}\big|\phi(\tfrac{t}{\tau})
                   -\phi(\tfrac{s}{\tau})\big||\psi(\tau)|d\tau &\mbox{if } t<s.
    \end{cases}.}

Consider first the case $s\in]t-\delta,t[$. The second inequality in
\eq{ln} implies that the measure of the interval $T:=]\frac{\ln t}{\ln s},1[$
is smaller than $\rho$. Thus, inequality \eq{T} holds with this set $T$. Therefore, for the first
term on the right hand side of \eq{II+}, using the estimates \eq{pq}, substituting $\tau=s^\sigma$,
and using the first inequality in \eq{ln}, we get
\Eq{B0}{
\int_s^t\tfrac{1}{\tau}\big|\phi(\tfrac{s}{\tau})\big||\psi(\tau)|d\tau
    &\leq\/\phi\/_p\/\psi\/_q\int_s^t\tfrac{1}{\tau}(\ln\tau-\ln s)^{p-1}(-\ln\tau)^{q-1}d\tau\\
    &\leq\/\phi\/_p\/\psi\/_q(-\ln s)^{p+q-1}
                \int_{\frac{\ln t}{\ln s}}^1(1-\sigma)^{p-1}\sigma^{q-1}d\sigma\\
    &<\frac{\/\phi\/_p\/\psi\/_q(-\ln s)^{p+q-1}\eps}{6\/\phi\/_p\/\psi\/_q(-\ln t)^{p+q-1}+1}
    <\frac{\eps}{6}.
}

To obtain an estimate for the second term on the right hand side
of \eq{II+} (when $s<t$), we use $\phi(x)=f(x)$ for $x\in H$ and
obtain
\Eq{EEE}{
 \int_{t}^1\tfrac{1}{\tau}&\big|\phi(\tfrac{t}{\tau})
                   -\phi(\tfrac{s}{\tau})\big||\psi(\tau)|d\tau
 \leq
 \int_{t}^1\tfrac{1}{\tau}
       \Big(\big|\phi(\tfrac{t}{\tau})-f(\tfrac{t}{\tau})\big|
       +\big|f(\tfrac{t}{\tau})-f(\tfrac{s}{\tau})\big|
       +\big|f(\tfrac{s}{\tau})-\phi(\tfrac{s}{\tau})\big|\Big)
  |\psi(\tau)|d\tau\\
  &\leq \int\limits_{]t,1[\setminus tH^{-1}}\hspace{-6mm}
          \tfrac1\tau \big|\phi(\tfrac{t}{\tau})-f(\tfrac{t}{\tau})\big||\psi(\tau)|d\tau
       + \int\limits_{]t,1[}\hspace{-2mm}
          \tfrac{1}{\tau}\big|f(\tfrac{t}{\tau})-f(\tfrac{s}{\tau})\big||\psi(\tau)|d\tau
       + \int\limits_{]t,1[\setminus sH^{-1}}\hspace{-6mm}
         \tfrac{1}{\tau}\big|f(\tfrac{s}{\tau})-\phi(\tfrac{s}{\tau})\big||\psi(\tau)|d\tau.
}
The first inequality in \eq{w} and the second estimate in \eq{ln} imply that, for $s\in]t-\delta,t]$,
\Eq{HHH}{
 \meas\big(]0,1[\setminus (1+(\ln s)^{-1}\widetilde{H})\big)
=(-\ln s)^{-1}\meas(]0,-\ln s[\setminus \widetilde{H})
<\frac{\rho\ln t}{2\ln s}<\frac{\rho(1+\rho)}{2}<\rho.
}
Thus \eq{T} holds with $T:=]0,1[\setminus (1+(\ln s)^{-1}\widetilde{H})$.
Using \eq{pq} and \eq{fp}, then substituting
$\tau=s^{\sigma}$ and finally applying inequality \eq{T}, we get
\Eq{B1}{
   \int\limits_{]t,1[\setminus sH^{-1}}\hspace{-6mm}
         \tfrac{1}{\tau}\big|f(\tfrac{s}{\tau})-\phi(\tfrac{s}{\tau})\big||\psi(\tau)|d\tau
   &\leq \int\limits_{]s,1[\setminus sH^{-1}}\hspace{-6mm}
         \tfrac{1}{\tau}\big(\big|f(\tfrac{s}{\tau})\big|+\big|\phi(\tfrac{s}{\tau})\big|\big)|\psi(\tau)|d\tau\\
   &\leq 2\/\phi\/_p\/\psi\/_q \int\limits_{]s,1[\setminus sH^{-1}}\hspace{-6mm}
           \tfrac1\tau (\ln\tau-\ln s)^{p-1}(-\ln\tau)^{q-1}d\tau\\
   &\leq 2\/\phi\/_p\/\psi\/_q (-\ln s)^{p+q-1}
         \int\limits_{]0,1[\setminus (1+(\ln s)^{-1}\widetilde{H})}\hspace{-6mm}
         (1-\sigma)^{p-1}\sigma^{q-1}d\sigma\\
   &\leq \frac{2\/\phi\/_p\/\psi\/_q (-\ln s)^{p+q-1}\eps}{6\/\phi\/_p\/\psi\/_q (-\ln t)^{p+q-1}+1}
      <\frac{\eps}{3}.
}
Applying this inequality for $s=t$, we also get
\Eq{B2}{
\int\limits_{]t,1[\setminus tH^{-1}}\hspace{-6mm}
         \tfrac{1}{\tau}\big|f(\tfrac{s}{\tau})-\phi(\tfrac{s}{\tau})\big||\psi(\tau)|d\tau
    <\frac{\eps}{3}.
}
Consider the second expression on the right hand side of
\eq{EEE}. We prove that
\Eq{NL}{
  \lim_{s\to t-0}\int\limits_{]t,1[}\hspace{-2mm}
   \tfrac{1}{\tau}\big|f(\tfrac{t}{\tau})-f(\tfrac{s}{\tau})\big||\psi(\tau)|d\tau=0.
}
By the continuity of $f$, the integrand pointwise converges to zero hence,
in view of Lebesgue's dominated convergence theorem, it suffices to show that
the integrand admits an integrable majorant which is independent of $s\in]t-\delta,t[$.

Using the inequality \eq{pq} and \eq{fp}, we get that, for all $\tau\in[t,1[$ and $s\in]t-\delta,t[$,
\Eq{*}
{
\tfrac{1}{\tau}\big|f(\tfrac{t}{\tau})-f(\tfrac{s}{\tau})\big||\psi(\tau)|
    &\leq \tfrac{1}{\tau}\Big(\big|f(\tfrac{t}{\tau})\big|+\big|f(\tfrac{s}{\tau})\big|\Big)|\psi(\tau)|\\
    &\leq \frac{\/\phi\/_p\/\psi\/_q}{\tau}
           \Big((\ln\tau-\ln t)^{p-1}+(\ln\tau-\ln s)^{p-1}\Big)(-\ln\tau)^{q-1}.
}
Now there are two cases. If $p\leq 1$ we have that the function $s\mapsto(\ln\tau-\ln s)^{p-1}$
is nondecreasing, hence
\Eq{*}
{
(\ln\tau-\ln s)^{p-1}\leq (\ln\tau-\ln t)^{p-1}\qquad (\tau\in[t,1[,\, s\in]0,t]).
}
This means that, in this case,
\Eq{*}
{
\tfrac{1}{\tau}\big|f(\tfrac{t}{\tau})-f(\tfrac{s}{\tau})\big||\psi(\tau)|
        \leq \frac{2\/\phi\/_p\/\psi\/_q}{\tau}(\ln\tau-\ln t)^{p-1}(-\ln\tau)^{q-1}
  \qquad (\tau\in[t,1[,\, s\in]0,t]).
}
Moreover, the right hand side is integrable with respect to $\tau$ because, with the substitution
$\tau=t^\sigma$, it follows that
\Eq{BB}{
   \int_t^1\tfrac{1}{\tau}(\ln\tau-\ln t)^{p-1}(-\ln\tau)^{q-1}d\tau
   =(-\ln t)^{p+q-1}B(p,q)<\infty.
}

When $p>1$, then the function $s\mapsto(\ln\tau-\ln s)^{p-1}$
is decreasing, hence
\Eq{*}
{
(\ln\tau-\ln s)^{p-1}\leq (\ln\tau-\ln(t-\delta))^{p-1}\leq (-\ln(t-\delta))^{p-1}
                      \qquad (\tau\in[t,1[,\,s\in]t-\delta,t]).
}
Thus, in this case, for all $\tau\in[t,1[,\,s\in]t-\delta,t]$,
\Eq{*}
{
\tfrac{1}{\tau}\big|f(\tfrac{t}{\tau})-f(\tfrac{s}{\tau})\big||\psi(\tau)|
             \leq \frac{\/\phi\/_p\/\psi\/_q}{\tau}\Big((\ln\tau-\ln t)^{p-1}
                   +(-\ln(t-\delta))^{p-1}\Big)(-\ln\tau)^{q-1}.
}
Again, the majorant is integrable because \eq{BB} holds, and
(substituting $\tau=t^{\sigma}$)
\Eq{*}
{
\int_{t}^1\tfrac{1}{\tau}(-\ln(t-\delta))^{p-1}(-\ln\tau)^{q-1}d\tau
        &=(-\ln(t-\delta))^{p-1}\int_{t}^1\tfrac{1}{\tau}(-\ln\tau)^{q-1}d\tau\\
        &=(-\ln(t-\delta))^{p-1}(-\ln t)^{q}B(1,q)<\infty.
}
Therefore, Lebesgue's Theorem can be applied and hence \eq{NL} holds.
Thus there exists $\delta^*\in]0,\delta]$, such that, for all $s\in]t-\delta^*,t[$,
\Eq{B3}
{
\int\limits_{]t,1[}\hspace{-2mm}
   \tfrac{1}{\tau}\big|f(\tfrac{t}{\tau})-f(\tfrac{s}{\tau})\big||\psi(\tau)|d\tau<\frac{\eps}{6}.
}
Combining the inequalities \eq{II+}, \eq{B0}, \eq{EEE}, \eq{B1}, \eq{B2}, and \eq{B3}, we get
\Eq{*}{
|(\phi*\psi)(t)&-(\phi*\psi)(s)|<\eps\qquad(s\in]t-\delta^*,t[),
}
which proves the left-continuity of $\phi*\psi$ at $t$.

To prove the right-continuity of $\phi*\psi$ at $t$, we apply \eq{II+} for $s\in]t,t+\delta[$.
The second inequality in \eq{ln} implies that
\Eq{*}
{
\meas\big(]\tfrac{\ln s}{\ln t},1[\big)\leq\frac{\rho}{1+\rho}<\rho.
}
Thus, inequality \eq{T} holds with the interval $T:=]\tfrac{\ln s}{\ln t},1[$.
Therefore, for the first term on the right hand side of \eq{II+} (using the estimates \eq{pq},
and substituting $\tau=t^\sigma$ in the evaluation of the integral), we get
\Eq{I1}{
\int_t^s\tfrac{1}{\tau}\big|\phi(\tfrac{t}{\tau})\big||\psi(\tau)|d\tau<\frac{\eps}{6}.
}
To obtain an estimate for the second term on the right hand side
of \eq{II+} (when $t<s$), we use $\phi(x)=f(x)$ for $x\in H$ and
obtain
\Eq{EEE+}{
 \int_{s}^1&\tfrac{1}{\tau}\big|\phi\big(\tfrac{t}{\tau}\big)
                   -\phi\big(\tfrac{s}{\tau}\big)\big||\psi(\tau)|d\tau
 \leq
 \int_{s}^1\tfrac{1}{\tau}
       \Big(\big|\phi\big(\tfrac{t}{\tau}\big)-f\big(\tfrac{t}{\tau}\big)\big|
       +\big|f\big(\tfrac{t}{\tau}\big)-f\big(\tfrac{s}{\tau}\big)\big|
       +\big|f\big(\tfrac{s}{\tau}\big)-\phi\big(\tfrac{s}{\tau}\big)\big|\Big)
  |\psi(\tau)|d\tau\\
  &\leq \hspace{-2mm}\int\limits_{]s,1[\setminus tH^{-1}}\hspace{-6mm}
          \tfrac1\tau \big|\phi\big(\tfrac{t}{\tau}\big)-f\big(\tfrac{t}{\tau}\big)\big||\psi(\tau)|d\tau
       + \int\limits_{]s,1[}\hspace{-2mm}
          \tfrac{1}{\tau}\big|f\big(\tfrac{t}{\tau}\big)-f\big(\tfrac{s}{\tau}\big)\big||\psi(\tau)|d\tau
       + \hspace{-2mm}\int\limits_{]s,1[\setminus sH^{-1}}\hspace{-6mm}
         \tfrac{1}{\tau}\big|f\big(\tfrac{s}{\tau}\big)-\phi\big(\tfrac{s}{\tau}\big)\big||\psi(\tau)|d\tau.
}
Applying an analogous argument as before, for $s\in]t,t+\delta[$, we can obtain the estimates
\Eq{I2}{
\int\limits_{]s,1[\setminus sH^{-1}}\hspace{-6mm}
         \tfrac{1}{\tau}\big|f\big(\tfrac{s}{\tau}\big)-\phi\big(\tfrac{s}{\tau}\big)\big||\psi(\tau)|d\tau
    <\frac{\eps}{3}\qquad\mbox{and}\qquad
\int\limits_{]s,1[\setminus tH^{-1}}\hspace{-6mm}
         \tfrac{1}{\tau}\big|f\big(\tfrac{s}{\tau}\big)-\phi\big(\tfrac{s}{\tau}\big)\big||\psi(\tau)|d\tau
    <\frac{\eps}{3}.
}

Consider the second expression on the right hand side of \eq{EEE+}.
We will prove that
\Eq{NL+}{
  \lim_{s\to t+0}\int\limits_{]s,1[}\hspace{-2mm}
   \tfrac{1}{\tau}\big|f\big(\tfrac{t}{\tau}\big)-f\big(\tfrac{s}{\tau}\big)\big||\psi(\tau)|d\tau=0.
}
First, with the substitution $\tau=\sigma^\frac{\ln s}{\ln t}$, for $s\in]t,t+\delta[$, we can obtain
\Eq{NL-}{
  \int\limits_{]s,1[}\hspace{-2mm}
   \tfrac{1}{\tau}\big|f\big(\tfrac{t}{\tau}\big)-f\big(\tfrac{s}{\tau}\big)\big||\psi(\tau)|d\tau
  &=\frac{\ln s}{\ln t}\int\limits_{]t,1[}\hspace{-2mm}
   \tfrac{1}{\sigma}\big|f\big(t\sigma^{-\frac{\ln s}{\ln t}}\big)-f\big(s\sigma^{-\frac{\ln s}{\ln t}}\big)\big|
        \big|\psi\big(\sigma^\frac{\ln s}{\ln t}\big)\big|d\sigma\\
 &\leq\int\limits_{]t,1[}\hspace{-2mm}
   \tfrac{1}{\sigma}\big|f\big(t\sigma^{-\frac{\ln s}{\ln t}}\big)-f\big(s\sigma^{-\frac{\ln s}{\ln t}}\big)\big|
        \big|\psi\big(\sigma^\frac{\ln s}{\ln t}\big)\big|d\sigma.
}
By the continuity of $f$ and the local boundedness of $\psi$ (which is a consequence of the
inequality \eq{pq}), the integrand on the right hand side of \eq{NL-} pointwise converges
to zero as $s\to t+0$, hence, in view of Lebesgue's dominated convergence theorem,
it suffices to show that
the integrand admits an integrable majorant which is independent of $s\in]t,t+\delta[$.
Using the inequality \eq{pq} and \eq{fp}, we get that, for all $\tau\in[t,1[$ and $s\in]t,t+\delta[$,
\Eq{*}
{
\tfrac{1}{\sigma}\big|f\big(t\sigma^{-\frac{\ln s}{\ln t}}\big)&
    -f\big(s\sigma^{-\frac{\ln s}{\ln t}}\big)\big|
        \big|\psi\big(\sigma^\frac{\ln s}{\ln t}\big)\big|\\
        &\leq \/\phi\/_p\/\psi\/_q\,\tfrac{1}{\sigma}
         \Big(\big(-\ln\big(t\sigma^{-\frac{\ln s}{\ln t}}\big)\big)^{p-1}
               +\big(-\ln\big(s\sigma^{-\frac{\ln s}{\ln t}}\big)\big)^{p-1}\Big)
                  \big(-\ln\big(\sigma^\frac{\ln s}{\ln t}\big)\big)^{q-1}.
}
Since $s\mapsto(-\ln s)^{q-1}$ is monotone on $]0,1[$, therefore, for $s\in]t,t+\delta[$,
we have
\Eq{*}
{
\big(-\ln\big(\sigma^\frac{\ln s}{\ln t}\big)\big)^{q-1}
 =\big(\tfrac{-\ln\sigma\ln s}{\ln t}\big)^{q-1}
\leq \max\big\{1,\big(\tfrac{\ln (t+\delta)}{\ln t}\big)^{q-1}\big\}(-\ln\sigma)^{q-1}.
}
Similarly, for $s\in]t,t+\delta[$ and $\sigma\in]t,1[$, we get
\Eq{*}{
\big(-\ln\big(s\sigma^{-\frac{\ln s}{\ln t}}\big)\big)^{p-1}
  =\big(\tfrac{\ln s}{\ln t}\big)^{p-1}(-\ln t+\ln\sigma)^{p-1}
  \leq \begin{cases}
        \big(\tfrac{\ln (t+\delta)}{\ln t}\big)^{p-1}(-\ln t+\ln\sigma)^{p-1}  &\mbox{if }p\leq 1,\\[2mm]
        (-\ln t+\ln\sigma)^{p-1}  &\mbox{if }p>1,
       \end{cases}
}
and
\Eq{*}{
\big(-\ln\big(t\sigma^{-\frac{\ln s}{\ln t}}\big)\big)^{p-1}
  =\big(-\ln t+\tfrac{\ln s\ln\sigma}{\ln t}\big)^{p-1}
  \leq\begin{cases}
       (-\ln t+\ln\sigma)^{p-1}&\mbox{if }p\leq 1,\\[2mm]
       (-\ln t)^{p-1}&\mbox{if }p>1.
       \end{cases}
}

Combining these inequalities, for $s\in]t,t+\delta[$ and $\sigma\in]t,1[$, we obtain
\Eq{*}{
\tfrac{1}{\sigma}&\big|f\big(t\sigma^{-\frac{\ln s}{\ln t}}\big)
    -f\big(s\sigma^{-\frac{\ln s}{\ln t}}\big)\big|
        \big|\psi\big(\sigma^\frac{\ln s}{\ln t}\big)\big|\\
  &\leq\begin{cases}
       \/\phi\/_p\/\psi\/_q
        \big(1+\big(\tfrac{\ln (t+\delta)}{\ln t}\big)^{p-1}\big)
        \max\big\{1,\big(\tfrac{\ln (t+\delta)}{\ln t}\big)^{q-1}\big\}
        \,\tfrac{1}{\sigma}(-\ln t+\ln\sigma)^{p-1}(-\ln\sigma)^{q-1}
           &\mbox{if }p\leq 1,\\[2mm]
       \/\phi\/_p\/\psi\/_q\max\big\{1,\big(\tfrac{\ln (t+\delta)}{\ln t}\big)^{q-1}\big\}
        \,\tfrac{1}{\sigma}\big((-\ln t+\ln\sigma)^{p-1}+(-\ln t)^{p-1}\big)(-\ln\sigma)^{q-1}
           &\mbox{if }p>1.
       \end{cases}
}
It is easy to check (by substituting $\sigma=t^\tau$) that the function on the right hand
side of this inequality is integrable with respect to $\sigma$ over $]t,1[$.
Therefore, Lebesgue's Theorem can be applied and hence \eq{NL+} holds.
Thus there exists $\delta^{**}\in]0,\delta^*]$, such that, for all $s\in]t,t+\delta^{**}[$,
\Eq{I3}
{
\int\limits_{]s,1[}\hspace{-2mm}
   \tfrac{1}{\tau}\big|f(\tfrac{t}{\tau})-f(\tfrac{s}{\tau})\big||\psi(\tau)|d\tau<\frac{\eps}{6}.
}
By summing up the respective sides of the inequalities \eq{I1}, \eq{EEE+}, \eq{I2}, and \eq{I3},
for all $s\in]t,t+\delta^{**}[$, we get
\Eq{*}{
|(\phi*\psi)(t)&-(\phi*\psi)(s)|<\eps,
}
which completes the proof of the right-continuity of $\phi*\psi$ at $t$.

To prove \eq{p+q}, let $t\in]0,1]$ be fixed.
Using \eq{pq} and substituting $\tau=t^s$, we get
\Eq{*} {
|(\phi*\psi)(t)|&\leq\int_t^1\tfrac{1}{\tau}\big|\phi(\tfrac{t}{\tau})\big||\psi(\tau)|d\tau
\leq \/\phi\/_p\/\psi\/_q|
     \int_t^1\tfrac{1}{\tau}\big(-\ln(\tfrac{t}{\tau})\big)^{p-1}(-\ln\tau)^{q-1}d\tau\\
&=\/\phi\/_p\/\psi\/_q(-\ln t)^{p+q-1}\int_0^1(1-s)^{p-1}s^{q-1}ds\\
&=B(p,q)\/\phi\/_p\/\psi\/_q(-\ln t)^{p+q-1}
 =\frac{\Gamma(p)\Gamma(q)}{\Gamma(p+q)}\/\phi\/_p\/\psi\/_q(-\ln t)^{p+q-1},
}
which proves the inclusion $\phi*\psi\in\Phi_{p+q}$ and the
inequality \eq{p+q}.

In the proof of \eq{prod} first we use Fubini's theorem and then
the variable $t/\tau$ is replaced by $s$:
 \Eq{*}{
\int_0^{1}(\phi*\psi)(t)dt&
  =\int_{0}^1\int_{t}^1\tfrac{1}{\tau}\phi(\tfrac{t}{\tau})\psi(\tau)d\tau dt
  =\int_{0}^1\int_{0}^{\tau}\tfrac{1}{\tau}\phi(\tfrac{t}{\tau})\psi(\tau)dt d\tau\\
  &=\int_{0}^1\psi(\tau)\Big(\tfrac{1}{\tau}\int_{0}^{\tau}\phi(\tfrac{t}{\tau})dt\Big)d\tau
  =\int_{0}^1\psi(\tau)d\tau\int_{0}^1\phi(s)ds.
}
\end{proof}

\Lem{Gamma}{Let $p>0$ be arbitrarily fixed, then, for all $x\in\R$,
\Eq{Gamma1}
{
\lim_{n\to\infty}\frac{x^n}{\Gamma(np)}=0
}
and the convergence is uniform on every compact interval of $\R$.}

\begin{proof}
To prove the lemma, we will show that the series $\sum \frac{x^n}{\Gamma(np)}$ is convergent
on $\R$. Using Cauchy's root test on this series and the Stirling formula for the $\Gamma$
function, we get that
\Eq{*}{
\lim_{n\to\infty}\sqrt[n]{ \frac{|x|^n}{\Gamma(np)}}
         =\lim_{n\to\infty} \frac{|x|}{\sqrt[n]{\Gamma(np)}}
         =|x|\lim_{n\to\infty} \frac{1}{\sqrt[n]{\sqrt{\tfrac{2\pi}{np}}(\tfrac{np}{e})^{np}}}
         =|x|\lim_{n\to\infty}\frac{e^p\sqrt[2n]{np}}{\sqrt[2n]{2\pi}(np)^p}=0.
}
This means that the series is absolute convergent on $\R$ and hence the convergence is uniform on
compact subsets of $\R$, which yields the statement.
\end{proof}

\Prp{p}{Let $p>0$ and $\phi\in\Phi_p$.
Define the sequence $\phi_n:]0,1[\to\R$ by the recursion
\Eq{phin} {
  \phi_1:=\phi,\qquad \phi_{n+1}:=\phi*\phi_{n}\qquad(n\in\N).
}
Then, for all $n\in\N$,
\Eq{int}{
 \phi_n\in\Phi_{np},\qquad
 \/\phi_n\/_{np}\leq \frac{(\Gamma(p))^n}{\Gamma(np)}\/\phi\/_{p}^n,\qquad
 \int_{0}^1\phi_{n}=\bigg(\int_{0}^1\phi\bigg)^n
}
and, for all $s\in]0,1[$,
\Eq{0} {
  \lim_{n\to\infty}\phi_n(s)=0
}
furthermore, for all $\delta\in]0,1[$, the convergence is uniform on $[\delta,1[$.}

\begin{proof} By the definition of $\phi_1$, \eq{int} holds trivially for $n=1$.
Assume that \eq{int} is valid for some $n$.
By \prp{pq} and the inductive assumption, $\phi_{n+1}=\phi*\phi_{n}\in\Phi_{p+np}$,
\Eq{*}{
 \/\phi_{n+1}\/_{(n+1)p}
  \leq \frac{\Gamma(p)\Gamma(np)}{\Gamma((n+1)p)}\/\phi\/_{p}\/\phi_{n}\/_{np}
  \leq \frac{\Gamma(p)\Gamma(np)}{\Gamma((n+1)p)}\/\phi\/_{p}
         \frac{(\Gamma(p))^n}{\Gamma(np)}\/\phi\/_{p}^n
  =\frac{(\Gamma(p))^{n+1}}{\Gamma((n+1)p)}\/\phi\/_{p}^{n+1},
}
and
\Eq{*}{
\int_{0}^1\phi_{n+1}&=\int_{0}^1\phi*\phi_n
                                 =\int_{0}^1\phi\int_{0}^1\phi_n
                                 =\int_{0}^1\phi\Big(\int_{0}^1\phi\Big)^{n}
                                 = \Big(\int_{0}^1\phi\Big)^{n+1},
}
which proves \eq{int} for $n+1$.
To prove \eq{0}, choose $n_0$ such that $n_0p\geq2$ and let $\delta\in]0,1[$.
Since the Gamma function is increasing on the interval $[2,\infty[$,
by \eq{int}, for all $n\geq n_0$ and $s\in[\delta,1[$, we have
\Eq{K}{
|\phi_{n+n_0}(s)|&\leq \/\phi_{n+n_0}\/_{(n+n_0)p}(-\ln s)^{(n+n_0)p-1}
\leq \frac{(\Gamma(p))^{n+n_0}}{\Gamma((n+n_0)p)}\/\phi\/_{p}^{n+n_0}(-\ln s)^{(n+n_0)p-1}\\
&=(-\ln s)^{n_0p-1}\Gamma(p)^{n_0}\/\phi\/_{p}^{n_0}
       \frac{\big(\Gamma(p)(-\ln s)^p\/\phi\/_{p}\big)^n}{\Gamma((n+n_0)p)}\\
&\leq(-\ln s)^{n_0p-1}\Gamma(p)^{n_0}\/\phi\/_{p}^{n_0}
       \frac{\big(\Gamma(p)(-\ln s)^p\/\phi\/_{p}\big)^n}{\Gamma(np)}\\
&\leq(-\ln \delta)^{n_0p-1}\Gamma(p)^{n_0}\/\phi\/_{p}^{n_0}
       \frac{\big(\Gamma(p)(-\ln \delta)^p\/\phi\/_{p}\big)^n}{\Gamma(np)}.
}
Letting $n\to\infty$ in \eq{K} and using \lem{Gamma}, we get \eq{0}. The estimate \eq{K} also
ensures the uniformity of the convergence on the compact subsets of $]0,1]$.
\end{proof}

\Prp{gphi}{Let $g:[0,1]\to\R$ be an upper bounded measurable function which is upper
semicontinuous at $0$. Let $p>0$, $\phi\in\Phi_p$ be a nonnegative function with
$\int_0^1\phi=1$ and define the the sequence $\phi_n:]0,1[\to\R$ by \eq{phin}. Then
\Eq{gphi}
{
\limsup_{n\to\infty}\int_0^1g(s)\phi_n(s)ds\leq g(0).
}}

\begin{proof}
Let $\eps>0$. If $g$ is upper semicontinuous at $0$, then there exists $\delta\in]0,1[$, such that
\Eq{g0}
{
g(s)<g(0)+\frac\eps2 \qquad \mbox{for all}\qquad s\in]0,\delta[,
}
and, by the upper boundedness, there exists $K>\max(0,-g(0))$, such that
\Eq{gK}
{
  g(s)\leq K\qquad (s\in [0,1]).
}
By \prp{p}, $\lim_{n\to\infty}\phi_n=0$ and the convergence is uniform also on
$[\delta,1[$. Hence, there exists $n_0\in\N$ such that, for all $n\geq n_0$,
\Eq{eps}
{
\phi_n(s)\leq \frac{\eps}{4K}\qquad (s\in[\delta,1[).
}
Using the third expression \eq{int} and $\int_0^1\phi=1$, we have that
\Eq{1}
{
\int_0^1\phi_n=\Big(\int_0^1\phi\Big)^n=1.
}
Applying \eq{1} and the nonnegativity of $\phi_n$, we obtain
\Eq{g1}
{
\int_0^1g(s)\phi_n(s)ds-g(0)= \int_0^\delta(g(s)-g(0))\phi_n(s)ds
           +\int_\delta^1 (g(s)-g(0))\phi_n(s)ds
}
To obtain an estimate for the first term of the right hand side, we use \eq{g0}
and \eq{1}. Then, for all $n\in\N$,
\Eq{g2}{
     \int_0^\delta(g(s)-g(0))\phi_n(s)ds
      \leq \int_0^\delta\frac\eps2\phi_n(s)ds
      \leq \frac\eps2\int_0^1\phi_n(s)ds=\frac\eps2.
}
Consider the second expression on the right hand side of \eq{g1}.
Then using \eq{gK} and \eq{eps}, we obtain, for all $n\geq n_0$,
\Eq{g3}{
\int_\delta^1(g(s)-g(0))\phi_n(s)ds
       \leq \int_\delta^1 2K\phi_n(s)ds
       \leq \int_\delta^1 2K\frac{\eps}{4K}ds=\frac\eps2(1-\delta)<\frac\eps2.
}
Combining the inequalities \eq{g1}, \eq{g2} and \eq{g3}, we get
that, for all $n\geq n_0$,
\Eq{*}
{
\int_0^1g(s)\phi_n(s)ds-g(0)<\eps,
}
which proves the statement.
\end{proof}

The proof of \thm{2} is based on a sequence of lemmas.

\Lem{1}{Let $\alpha_H:D^*\to\R$ be even, $\rho:[0,1]\to\R$
be integrable and $\lambda\in]0,1[$. Then every $f:D\to\R$ lower hemicontinuous function
satisfying the approximate Hermite--Hadamard inequality \eq{H}, fulfills
\Eq{L1+}{
  \int_0^1 \Big(f\big(\tfrac{1+s}{2}x+\tfrac{1-s}{2}y\big)
      &+f\big(\tfrac{1-s}{2}x+\tfrac{1+s}{2}y\big)\Big)
        \frac{\rho\big(\tfrac{1+s}{2}\big)+\rho\big(\tfrac{1-s}{2}\big)}{2}ds\\
  &\leq f(x)+f(y)+2\alpha_H(x-y)\qquad(x,y\in D).
}}

\begin{proof}
Changing the role of $x$ and $y$ in \eq{H}, then adding the respective sides of the inequality so obtained
and the original inequality \eq{H}, by the evenness of $\alpha_H$, we get that
\Eq{i1}
{
 \int_0^1 \big(f(tx+(1-t)y)+f((1-t)x+ty)\big)\rho(t)dt
   \leq f(x)+f(y)+2\alpha_H(x-y)\quad(x,y\in D).
}
Replacing $t$ by $1-t$ in the integral on the left hand side of \eq{i1}, it follows that
\Eq{i2}
{
\int_0^1 \big(f((1-t)x+ty)+f(tx+(1-t)y)\big)\rho(1-t)dt
   \leq f(x)+f(y)+2\alpha_H(x-y)\quad(x,y\in D),
}
hence, adding the respective sides of the inequalities \eq{i1} and \eq{i2}, we obtain
\Eq{i12}
{
 \int_0^1 \big(f(tx+(1-t)y)+f((1-t)x+ty)\big)\tfrac{\rho(t)+\rho(1-t)}{2}dt
   \leq f(x)+f(y)+2\alpha_H(x-y)\quad(x,y\in D).
}
Finally, substituting $t:=\tfrac{1+s}{2}$ in the integral on the left hand side of \eq{i12},
we arrive at
\Eq{i12h}{
  \frac12\int_{-1}^1 \Big(f\big(\tfrac{1+s}{2}x+\tfrac{1-s}{2}y\big)
      &+f\big(\tfrac{1-s}{2}x+\tfrac{1+s}{2}y\big)\Big)
        \frac{\rho\big(\tfrac{1+s}{2}\big)+\rho\big(\tfrac{1-s}{2}\big)}{2}ds\\
  &\leq f(x)+f(y)+2\alpha_H(x-y)\qquad(x,y\in D).
}
Since the integrand on the left hand side of \eq{i12h}
is even, this inequality reduces to \eq{L1+}, which completes the proof of the lemma.
\end{proof}

In what follows, we examine the Hermite--Hadamard inequality \eq{L1+}.

\Lem{2}{Let $\rho:[0,1]\to\R_+$ be integrable with $\int_0^1\rho=1$ and there
exist $c\geq0$ and $p>0$ such that \eq{c} holds. Define $\phi:]0,1[\to\R$ by
\Eq{L2+}
{
\phi(s):=\frac{\rho\big(\tfrac{1+s}{2}\big)+\rho\big(\tfrac{1-s}{2}\big)}{2} \qquad (s\in ]0,1[).
}
Then $\phi\in\Phi_p$ and $\int_0^1\phi=1.$}

\begin{proof} Inequality \eq{c} results that
\Eq{*}{
  \phi(s)=\frac{\rho\big(\tfrac{1+s}{2}\big)+\rho\big(\tfrac{1-s}{2}\big)}{2}
    \leq c(-\ln s)^{p-1}\qquad(s\in]0,1[),
}
which proves $\phi\in\Phi_p$. The equality $\int_0^1\phi=1$ is an automatic
consequence of the assumption $\int_0^1\rho=1$.
\end{proof}

\Lem{3}{Let $p,q>0$ and $\phi\in\Phi_p$, $\psi\in\Phi_q$ be nonnegative functions. Let
$\alpha:D^*\to\R$ and let $\beta:D^*\to\R$ be a radially upper semicontinuous function. Assume that
a lower hemicontinuous function $f:D\to\R$ satisfies the approximate Hermite--Hadamard inequalities
\Eq{A}{
\int_{0}^1 \Big(f\big(\tfrac{1+s}{2}x+\tfrac{1-s}{2}y\big)
      +f\big(\tfrac{1-s}{2}x+\tfrac{1+s}{2}y\big)\Big)
        \phi(s)ds
  \leq f(x)+f(y)+2\alpha(x-y)\quad(x,y\in D),
}
and
\Eq{B}{
\int_{0}^1 \Big(f\big(\tfrac{1+s}{2}x+\tfrac{1-s}{2}y\big)
      +f\big(\tfrac{1-s}{2}x+\tfrac{1+s}{2}y\big)\Big)
       \psi(s)ds
  \leq f(x)+f(y)+2\beta(x-y)\quad(x,y\in D).
}
Then $f$ also satisfies the inequality
\Eq{AB}{
\int_{0}^1 \Big(f\big(\tfrac{1+s}{2}x+\tfrac{1-s}{2}y\big)
      &+f\big(\tfrac{1-s}{2}x+\tfrac{1+s}{2}y\big)\Big)
        (\phi*\psi)(s)ds\\
  &\leq f(x)+f(y)+2\alpha(x-y)+2\int_0^1\beta(t(x-y))\phi(t)dt\quad(x,y\in D).
}
}

\begin{proof}
Assume that $f:D\to\R$ satisfies the inequalities \eq{A} and \eq{B}.
To prove \eq{AB}, let $x,y\in D$ be fixed. Applying \eq{B} for the elements
$\frac{1+t}{2}x+\frac{1-t}{2}y,\frac{1-t}{2}x+\frac{1+t}{2}y\in D$, we obtain
\Eq{*}
{\int_{0}^1 \Big(f\big(\tfrac{1+ts}{2}x&+\tfrac{1-ts}{2}y\big)
      +f\big(\tfrac{1-ts}{2}x+\tfrac{1+ts}{2}y\big)\Big)
       \psi(s)ds\\
  &\leq f\big(\tfrac{1+t}{2}x+\tfrac{1-t}{2}y\big)+f\big(\tfrac{1-t}{2}x+\tfrac{1+t}{2}y\big)
      +2\beta(t(x-y)).
}
Multiplying this inequality by $\phi(t)$ and integrating the functions on both sides with 
respect to $t$ on $]0,1[$, then using that $f$ also satisfies \eq{A}, we get that
\Eq{AA}{
\int_{0}^1\int_{0}^1&\Big(f\big(\tfrac{1+ts}{2}x+\tfrac{1-ts}{2}y\big)
      +f\big(\tfrac{1-ts}{2}x+\tfrac{1+ts}{2}y\big)\Big)
       \psi(s)ds\phi(t)dt\\
&\leq \int_{0}^1\Big(f\big(\tfrac{1+t}{2}x+\tfrac{1-t}{2}y\big)
     +f\big(\tfrac{1-t}{2}x+\tfrac{1+t}{2}y\big)\Big)\phi(t)dt
     +2\int_{0}^1\beta(t(x-y))\phi(t)dt\\
&\leq f(x)+f(y)+2\alpha(x-y)+2\int_{0}^1\beta(t(x-y))\phi(t)dt.
}
Now we compute the left hand side of the previous inequality.
Substituting $s=\frac{\tau}{t}$ and using also Fubini's theorem, we obtain
\Eq{=}{
\int_{0}^1\int_{0}^1\Big(f\big(\tfrac{1+ts}{2}x&+\tfrac{1-ts}{2}y\big)
      +f\big(\tfrac{1-ts}{2}x+\tfrac{1+ts}{2}y\big)\Big)
       \psi(s)ds\phi(t)dt\\
       &=\int_0^1\int_0^t\Big(f\big(\tfrac{1+\tau}{2}x+\tfrac{1-\tau}{2}y\big)
      +f\big(\tfrac{1-\tau}{2}x+\tfrac{1+\tau}{2}y\big)\Big)
          \psi(\tfrac{\tau}{t})\tfrac{1}{t}d\tau\phi(t)dt\\
       &=\int_0^1\Bigg(\Big(f\big(\tfrac{1+\tau}{2}x+\tfrac{1-\tau}{2}y\big)
      +f\big(\tfrac{1-\tau}{2}x+\tfrac{1+\tau}{2}y\big)\Big)
          \int_{\tau}^1\tfrac{1}{t}\psi(\tfrac{\tau}{t})\phi(t)dt\Bigg)d\tau\\
       &=\int_0^1\Big(f\big(\tfrac{1+\tau}{2}x+\tfrac{1-\tau}{2}y\big)
      +f\big(\tfrac{1-\tau}{2}x+\tfrac{1+\tau}{2}y\big)\Big)(\psi*\phi)(\tau)d\tau}.
Combining \eq{AA} and \eq{=}, the inequality \eq{AB} follows, which completes the proof.
\end{proof}

\Lem{3.5}{Let $\phi:]0,1[\to\R_+$ be an integrable function and let
$\beta:D^*\to\R$ be a radially upper semicontinuous. Then, the function
$\gamma:D^*\to\R$ defined by
\Eq{gam}{
  \gamma(u):=\int_0^1\beta(tu)\phi(t)dt\qquad(u\in D^*)
}
is also radially upper semicontinuous on $D^*$.}

\begin{proof} To prove that $\gamma$ defined by \eq{gam} is radially upper semicontinuous
at $u_0\in D^*$, let $s_n\to s_0$ be an arbitrary sequence in $[0,1]$.
We have that
\Eq{*}{
\beta(ts_nu_0)\leq\sup_{\tau\in [0,1]}\beta(\tau u_0)=:K\qquad(t\in[0,1],\,n\in\N),
}
thus, $K\phi$ is an integrable majorant for the sequence of functions
$t\mapsto\beta(ts_nu_0)\phi(t)$. Using Fatou's lemma and the radial upper semicontinuity of
$\beta$, we get that
\Eq{*}{
\limsup_{n\to\infty}\gamma(s_nu_0)
               &=\limsup_{n\to\infty}\int_0^1\beta(ts_nu_0)\phi(t)dt\\
               &\leq \int_0^1\limsup_{n\to\infty}\beta(ts_nu_0)\phi(t)dt
               \leq\int_0^1\beta(ts_0u_0)\phi(t)dt
               =\gamma(s_0u_0),
}
which proves the statement.
\end{proof}

\Lem{4}{Let $p>0$, $\phi\in\Phi_p$ be a nonnegative and
$\alpha_H:D^*\to\R$ be a radially upper semicontinuous function.
If $f:D\to\R$ is lower hemicontinuous and fulfills the approximate Hermite--Hadamard inequality
\Eq{H'}{
\int_{0}^1\!\! \Big(f\big(\tfrac{1+s}{2}x+\tfrac{1-s}{2}y\big)
      +f\big(\tfrac{1-s}{2}x+\tfrac{1+s}{2}y\big)\Big)\phi(s)ds
  \leq f(x)+f(y)+2\alpha_H(x-y)\quad(x,y\in D),
}
then, for all $n\in\N$, the function $f$ also satisfies the Hermite--Hadamard inequality
\Eq{H''}{
\int_{0}^1\!\!\! \Big(f\big(\tfrac{1+s}{2}x+\tfrac{1-s}{2}y\big)
      +f\big(\tfrac{1-s}{2}x+\tfrac{1+s}{2}y\big)\Big)\phi_n(s)ds
  \leq f(x)+f(y)+2\alpha_n(x-y)\quad(x,y\in D),
}
where the sequences $\phi_n:]0,1[\to\R_+$ and $\alpha_n:D^*\to\R_+$ are defined by
\eq{phin} and
\Eq{alpha}{
\alpha_1=\alpha_H, \qquad \alpha_{n+1}(u)
     =\int_0^1\alpha_n(tu)\phi(t)dt+\alpha_H(u)\qquad(u\in D^*),
}
respectively.}

\begin{proof} We note that, by \lem{3.5}, the sequence of functions $(\alpha_n)$
is well-defined and $\alpha_n$ is radially lower semicontinuous for all $n\in N$.

Let $x,y\in D$ and $p>0$. To prove \eq{H''}, we use induction on $n\in\N$.
For $n=1$, we have \eq{H'}. Assume that \eq{H''} holds for $n\in\N$.
Since $\phi\in\Phi_p$, by \prp{p}, we have that $\phi_n\in\Phi_{np}$.
The function $f$ satisfies \eq{H'} and also \eq{H''}, for $n\in\N$.
Thus, in \lem{3}, \eq{A} holds with the functions $\phi$ and $\alpha:=\alpha_1$.
Furthermore, by the inductive assumption, also in \lem{3},
\eq{B} holds with the functions $\psi:=\phi_n$ and $\beta:=\alpha_n$, for $n\in\N$.
Hence the function $f$ also fulfills the Hermite--Hadamard inequality \eq{AB}, which results,
\Eq{*}{
\int_{0}^1\Big(f\big(\tfrac{1+s}{2}x+\tfrac{1-s}{2}y\big)
      &+f\big(\tfrac{1-s}{2}x+\tfrac{1+s}{2}y\big)\Big)(\phi*\phi_n)(s)ds\\
  &\leq f(x)+f(y)+2\alpha_H(x-y)+2\int_{0}^1\alpha_n(t(x-y))\phi(t)dt,
}
which is the case $n+1$.
\end{proof}

\Lem{5}{Let $p>0$ and $\phi\in\Phi_p$ be a nonnegative function with $\int_0^1\phi(t)dt=1$
and $f:D\to\R$ be lower hemicontinuous. Then
\Eq{LIM1}{
\liminf_{n\to\infty}\int_0^1\Big(f\big(\tfrac{1+s}{2}x+\tfrac{1-s}{2}y\big)
      +f\big(\tfrac{1-s}{2}x+\tfrac{1+s}{2}y\big)\Big)\phi_n(s)ds
   \geq 2f\Big(\frac{x+y}{2}\Big)\qquad(x,y\in D).
}}

\begin{proof}
To prove the statement, let $x,y\in D$ and $p>0$.
Define $g_{x,y}:[0,1]\to\R$ by
\Eq{*}{
g_{x,y}(s):&=f\big(\tfrac{1+s}{2}x+\tfrac{1-s}{2}y\big)
               +f\big(\tfrac{1-s}{2}x+\tfrac{1+s}{2}y\big)\qquad (s\in[0,1]).
}
The lower hemicontinuity of $f$ implies that $g_{x,y}$ is lower semicontinuous and
hence lower bounded on $[0,1]$.
Thus, we can apply \prp{gphi}, for $g:=-g_{x,y}$ and $\phi\in\Phi_p$, which yields that
\Eq{*}{
\liminf_{n\to\infty}\int_0^1g_{x,y}(s)\phi_n(s)ds\geq g_{x,y}(0).
}
This inequality is equivalent to \eq{LIM1}.
\end{proof}

\Lem{6}{Let $p>0$ and $\phi\in\Phi_p$ be a nonnegative function,
and $\alpha_H:D^*\to\R$ be a radially upper semicontinuous function.
Then, for all $n\in\N$, the function $\alpha_n:D^*\to\R$ defined by \eq{alpha}
is radially upper semicontinuous and the sequence $(\alpha_n)$ is nondecreasing
[nonincreasing], whenever $\alpha_H$ is nonnegative [nonpositive].
Furthermore, if $\alpha_J:D^*\to\R$ is a radially lower semicontinuous
solution of the functional inequality
\Eq{aa}{
\alpha_J(u)\geq \int_0^1\alpha_J(su)\phi(s)ds+\alpha_H(u)\qquad (u\in D^*),
}
then
\Eq{limalpha}{
\limsup_{n\to\infty}\alpha_n(u)\leq \alpha_J(u)-\alpha_J(0)+\alpha_H(0)\qquad (u\in D^*).
}}

\begin{proof} The statement about the radial upper semicontinuity directly follows from \lem{3.5}.

Assume first that $\alpha_H$ is nonnegative. We will prove by induction on $n\in \N$, that the
sequence $(\alpha_n)$ is nondecreasing, i.e.,
\Eq{ia}
{
\alpha_{n+1}\geq\alpha_{n}\qquad (n\in\N).
}
For $n=1$, by the nonnegativity of $\alpha_1=\alpha_H$, we have that
\Eq{*}
{
\alpha_2(u)=\int_0^1\alpha_1(su)\phi(s)ds+\alpha_H(u)
   \geq \alpha_1(u)\qquad (u\in D^*).
}
Assume that \eq{ia} holds for some $n\in\N$ and consider the case $n+1$. Using the definition of
$\alpha_{n+1}$, the inductive assumption and the nonnegativity of $\alpha_n$, we get that
\Eq{*}
{
\alpha_{n+2}(u)
    =\int_{0}^1\alpha_{n+1}(su)\phi(s)ds+\alpha_H(u)
     \geq \int_{0}^1\alpha_{n}(su)\phi(s)ds+\alpha_H(u)
     =\alpha_{n+1}(u)
     \quad (u\in D^*).
}
Analogously, if $\alpha_H$ is nonpositive, we can obtain that
the sequence $(\alpha_n)$ is nonincreasing.

To prove \eq{limalpha}, let $\alpha_J:D^*\to\R$ be a radially lower semicontinuous solution
of \eq{aa}. Subtracting the respective sides of the inequalities \eq{aa} from \eq{alpha}, for the
sequence of functions $g_n:=\alpha_n-\alpha_J$, we obtain
\Eq{*}{
g_{n+1}(u)\leq \int_0^1g_n(su)\phi(s)ds\qquad (u\in D^*,\,n\in\N),
}
We obviously have that $g_n$ is also radially upper semicontinuous.

Iterating this inequality, similarly as in \lem{3} and \lem{4}, it can be proved that
\Eq{betaJ}
{
g_{n+1}(u)\leq \int_0^1g_1(su)\phi_n(s)ds\qquad (u\in D^*,\,n\in\N),
}
where $\phi_n$ is defined by \eq{phin}. Taking the limsup as $n\to\infty$ in \eq{betaJ},
by \prp{gphi}, we get that
\Eq{*}
{
\limsup_{n\to\infty}g_{n+1}(u)
    \leq \limsup_{n\to\infty}\int_0^1g_1(su)\phi_n(s)ds
    \leq g_1(0)=\alpha_H(0)-\alpha_J(0)\qquad (u\in D^*),
}
which immediately yields \eq{limalpha}.
\end{proof}

\begin{proof}[Proof of \thm{2}] Assume that the conditions of \thm{2} hold and $f:D\to\R$
is an upper semicontinuous solution of \eq{H}. Then by \lem{1}, $f$ also fulfills \eq{H'}, where
$\phi:]0,1[\to\R_+$ is defined by \eq{L2+}. Then, by \lem{2}, $\phi\in\Phi_p$ and $\int_0^1\phi=1.$
Using \lem{4}, we get that \eq{H''} also holds, where, for all $n\in\N$, $\alpha_{n}:D^*\to\R$ is
defined by \eq{alpha}. Since $\alpha_J$ satisfies the functional inequality \eq{a}, thus
applying Fubini's theorem, then substituting $s:=1-2t$ and $s:=2t-1$, we get that
\Eq{*}{
\alpha_J(u)&\geq \int_0^{1/2}\alpha_J((1-2t)u)\rho(t)dt
            +\int_{1/2}^1\alpha_J((2t-1)u)\rho(t)dt+\alpha_H(u)\\
             &=\int_0^{1}\alpha_J(su)\rho(\tfrac{1-s}{2})\tfrac12ds
                 +\int_0^{1}\alpha_J(su)\rho(\tfrac{1+s}{2})\tfrac12ds+\alpha_H(u)\\
             &=\int_0^{1}\alpha_J(su)\frac{\rho(\tfrac{1+s}{2})+\rho(\tfrac{1-s}{2})}{2}ds+\alpha_H(u)
             =\int_0^{1}\alpha_J(su)\phi(s)ds+\alpha_H(u),
}
which means that \eq{aa} also holds. Taking the liminf as $n\to\infty$ in \eq{H''}, using also
\lem{5}, \lem{6} and $\alpha_H(0)\leq\alpha_J(0)$, we get that
\Eq{*}{
  2f\Big(\frac{x+y}{2}\Big)&\leq f(x)+f(y)+2\liminf_{n\to\infty}\alpha_n(x-y)
  \leq f(x)+f(y)+2\limsup_{n\to\infty}\alpha_n(x-y)\\
 &\leq f(x)+f(y)+2\big(\alpha_J(x-y)+\alpha_H(0)-\alpha_J(0)\big)\leq f(x)+f(y)+2\alpha_J(x-y).
}
Hence \eq{J} holds, which completes the proof of \thm{2}.
\end{proof}

In what follows, we examine the case, when $X$ is a normed space and
$\alpha_H$ is a linear combination of the powers
of the norm with positive exponents, i.e., if $\alpha_H$ is of the form
\Eq{DHP}{
  \alpha_H(u):=\int\limits_{]0,\infty[} \|u\|^qd\mu_H(q)\qquad (u\in D^*),
}
where $\mu_H$ is a signed Borel measure on the interval $]0,\infty[$.
An important particular case is when $\mu_H$ is of the form $\sum_{i=1}^k c_i\delta_{q_i}$,
where $c_1,\dots,c_k\in\R$, $q_1,\dots,q_k>0$ and $\delta_q$
denotes the Dirac measure concentrated at $q$.

\Thm{A1}{Let $\rho:[0,1]\to\R$ be integrable with $\int_0^1\rho=1$ and assume that
there exist $c\geq0$ and $p>0$ such that \eq{c} holds.
Let $\lambda\in[0,1]$ and $\mu_H$ be a signed Borel measure on $]0,\infty[$ such that
\Eq{mh}{
  \int\limits_{]0,\infty[}\|u\|^qd|\mu_H|(q)<\infty
   \qquad(u\in D^*)
}
and
\Eq{mh+}{
  \int\limits_{]0,\infty[}\bigg(\int_0^1\big(1-|1-2t|^q\big)\rho(t)dt\bigg)^{-1}d|\mu_H|(q)<\infty.
}
Assume that $f:D\to\R$ is lower hemicontinuous and satisfies the
Hermite--Hadamard type inequality
\Eq{A1H}{
\int_0^1f(tx+(1-t)y)\rho(t)dt\leq\lambda f(x)+(1-\lambda)f(y)
    +\int\limits_{]0,\infty[} \|x-y\|^qd\mu_H(q)\qquad (x,y\in D).
}
Then $f$ also fulfils the Jensen type inequality
\Eq{A1J}{
f\Big(\frac{x+y}{2}\Big)\leq \frac{f(x)+f(y)}{2}
    +\!\!\int\limits_{]0,\infty[}\!\!\bigg(\int_0^1\big(1-|1-2t|^q\big)\rho(t)dt\bigg)^{-1}\|x-y\|^qd\mu_H(q)
     \quad (x,y\in D).
}}

\begin{proof} By \thm{2}, it suffices to show that the function
\Eq{*}{
  \alpha_J(u)
  :=\int\limits_{]0,\infty[}\!\!\bigg(\int_0^1\big(1-|1-2t|^q\big)\rho(t)dt\bigg)^{-1}\|u\|^qd\mu_H(q)
     \quad (u\in D^*)
}
is well-defined and satisfies \eq{a} with equality where $\alpha_H:D^*\to\R$ is defined by
\Eq{*}{
  \alpha_H(u)
  :=\int\limits_{]0,\infty[}\|u\|^qd\mu_H(q)   \qquad (u\in D^*)
}

To see that, for all $u\in D^*$, $\alpha_J(u)$ is finite, we distinguish two cases.
If $\|u\|\leq 1$, then $\|u\|^q\leq 1$ for all $q>0$, and hence, by assumption \eq{mh+},
\Eq{*}{
 |\alpha_J(u)|
 \leq \int\limits_{]0,\infty[}\!\!\bigg(\int_0^1\big(1-|1-2t|^q\big)\rho(t)dt\bigg)^{-1}d|\mu_H|(q)<\infty.
}
Now let $\|u\|>1$. Then, the functions $q\mapsto \|u\|^q$ and
$q\mapsto\int_0^1\big(1-|1-2t|^q\big)\rho(t)dt$ are increasing functions, hence
\Eq{*}{
 &|\alpha_J(u)|\\
  &\leq \|u\|\int\limits_{]0,1]}\!\!\bigg(\int_0^1\big(1-|1-2t|^q\big)\rho(t)dt\bigg)^{-1}d|\mu_H|(q)
   +\bigg(\int_0^1\big(1-|1-2t|\big)\rho(t)dt\bigg)^{-1}
        \!\!\!\!\int\limits_{]1,\infty[}\!\!\|u\|^qd|\mu_H|(q),
}
which is again finite by conditions \eq{mh} and \eq{mh+}.

To prove that $\alpha_J$ satisfies \eq{a}, using $\int_0^1\rho=1$, we compute
\Eq{*}{
&\int_0^1\alpha_J(|1-2s|u)\rho(s)ds+\alpha_H(u)\\
  &=\int_0^1\int\limits_{]0,\infty[}
     \!\!\bigg(\int_0^1\big(1-|1-2t|^q\big)\rho(t)dt\bigg)^{-1}\||1-2s|u\|^qd\mu_H(q)\rho(s)ds
   +\int\limits_{]0,\infty[} \|u\|^qd\mu_H(q)\\
  &=\!\!\int\limits_{]0,\infty[}\!\!\bigg(
      \frac{\int_0^1\!\!|1-2s|^q\rho(s)ds}{\int_0^1\big(1-|1-2t|^q\big)\rho(t)dt}+1\bigg)\|u\|^q
   d\mu_H(q)
  =\!\!\int\limits_{]0,\infty[}\!\!
      \frac{\|u\|^q}{\int_0^1\big(1-|1-2t|^q\big)\rho(t)dt}d\mu_H(q)=\alpha_J(u),
}
which proves that \eq{a} holds with equality.
\end{proof}

\Cor{A1}{
Let $\lambda\in[0,1]$, $a\in\R$ and $q>0$.
Assume that $f:D\to\R$ is lower hemicontinuous and satisfies the
Hermite--Hadamard type inequality
\Eq{*}
{
\int_0^1f(tx+(1-t)y)dt\leq\lambda f(x)+(1-\lambda)f(y)
    +a\|x-y\|^q\qquad (x,y\in D).
}
Then $f$ also fulfils the Jensen type inequality
\Eq{*}
{
f\Big(\frac{x+y}{2}\Big)\leq \frac{f(x)+f(y)}{2}
    +a\frac{q+1}{q} \|x-y\|^q  \qquad (x,y\in D).
}}

\begin{proof} Observe that the constant weight function $\rho\equiv 1$ satisfies the assumptions of
\thm{A1} with $c=p=1$. Also, with $\mu_H:=a\delta_q$, conditions \eq{mh} and \eq{mh+} hold
trivially. Thus, the conclusion of \thm{A1} is valid with
\Eq{*}{
\alpha_J(u)=\frac{a\|u\|^q}{\int_0^1\big(1-|1-2t|^q\big)dt}
    &=\frac{a\|u\|^q}{\int_0^{\frac12}\big(1-(1-2t)^q\big)dt+\int_{\frac12}^1\big(1-(2t-1)^q\big)dt}\\
    &=\frac{a\|u\|^q}{2\big(\tfrac12-\tfrac1{2(q+1)}\big)}
    =a\frac{q+1}{q}\|u\|^q\qquad(u\in D^*),
}
which proves the statement.
\end{proof}

\textbf{Acknowledgement.} The authors wish to thank the anonymous referee for the careful reading and
the useful advices.

%\bibliography{publ,funcequ,control}
%\bibliographystyle{plain}

\end{document}